\documentclass{article}

\usepackage{amsmath, amssymb, amsfonts, amsbsy, latexsym}
\usepackage{graphicx,stmaryrd}
\usepackage{epstopdf}
\usepackage{setspace, color}

\newtheorem{theorem}{Theorem}[section]

\newtheorem{lemma}[theorem]{Lemma}
\newtheorem{proposition}[theorem]{Proposition}

\newtheorem{remark}[theorem]{Remark}
\newtheorem{example}[theorem]{Example}
\newtheorem{definition}[theorem]{Definition}

\def\cB{\mathcal{B}}

\def\cF{\mathcal{F}}

\def\cR{\mathcal{R}}
\def\C{\mathbb{C}}
\def\R{\mathbb{R}}

\def\NN{\mathbb{N}}

\def\RR{\mathbb{R}}
\def\cC{\mathbb{C}}

\def\<{\langle}
\def\>{\rangle}

\def\epsilon{\varepsilon}

\def\l{\lambda}

\newcommand{\ip}[2]{\langle#1,#2\rangle}
\newcommand{\abs}[1]{\left| #1 \right|}
\newcommand{\absip}[2]{\left| \left\langle#1,#2\right\rangle \right|}
\newcommand{\norm}[1]{\left\lVert#1\right\rVert}

\newcommand{\fc}{{\mathcal F}}

\newcommand{\hc}{{H}}
\newcommand{\hch}{{\hat{\hc}}}

\newcommand{\phik}{{\varphi_k}}
\newcommand{\outp}[2]{{\llbracket #1,#2 \rrbracket }}
\newcommand{\ac}{\mathcal{A}}
\newcommand{\A}{\mathcal{A}}
\renewcommand{\SS}{\mathcal{S}}

\renewcommand{\j}{{\mathbf j}}
\newcommand{\Sym}{\textup{\mbox{Sym}}}
\newcommand{\trace}{\textup{\mbox{trace}}}
\newcommand{\Lip}{\textup{\mbox{Lip}}}
\newcommand{\real}{\textup{\mbox{real}}}
\newcommand{\imag}{\textup{\mbox{imag}}}

\title{On Lipschitz Analysis and Lipschitz Synthesis for the Phase Retrieval Problem}

\author{Radu Balan$^{(1)}$ , Dongmian Zou$^{(2)}$\\
$^{(1)}$ {\small Department of Mathematics and Center for Scientific Computation and Mathematical Modeling} \\
{\small University of Maryland, College Park, MD 20742, USA} \\
$^{(2)}$ {\small Applied Mathematics, Applied Statistics and Scientific Computing Program} \\
{\small University of Maryland, College Park, MD 20742, USA} \\
}

\begin{document}

\maketitle

\begin{abstract}
In this paper we prove two results regarding reconstruction from magnitudes of frame coefficients (the so called "phase retrieval problem"). First we show that phase retrievability as an algebraic property implies that nonlinear maps are bi-Lipschitz with respect to appropriate metrics on the quotient space. Second we prove that reconstruction can be
performed using Lipschitz continuous maps. Specifically we show that when nonlinear analysis maps $\alpha,\beta:\hch\rightarrow\R^m$ are injective, with $\alpha(x)=(|\ip{x}{f_k}|)_{k=1}^m$ and $\beta(x)=(|\ip{x}{f_k}|^2)_{k=1}^m$, 
where $\{f_1,\ldots,f_m\}$ is a frame for a Hilbert space $\hc$ and $\hch=\hc/T^1$, then $\alpha$ is bi-Lipschitz
with respect to the class of "natural metrics" $D_p(x,y)= min_{\varphi}\norm{x-e^{i\varphi}y}_p$, whereas $\beta$
is bi-Lipschitz with respect to the class of matrix-norm induced metrics $d_p(x,y)=\norm{xx^*-yy^*}_p$.  
Furthermore, there exist left inverse maps $\omega,\psi:\R^m\rightarrow \hch$ of $\alpha$ and $\beta$ respectively, 
that are Lipschitz continuous with respect to the appropriate metric. Additionally we obtain the
Lipschitz constants of these inverse maps in terms of the lower Lipschitz constants of $\alpha$ and $\beta$. 
Surprisingly the increase in Lipschitz constant is a relatively small factor, independent of
the space dimension or the frame redundancy.
\end{abstract}

\section{Introduction}\label{sec1}

Assume $\fc=\{f_1,f_2,\ldots,f_m\}$ is a frame (that is a spanning set) for the $n$-dimensional Hilbert space $H$. In this paper $H$ can be
a real or complex Hilbert space. The results in Section \ref{sec3} apply to both cases, and the constants have the same form.

On $H$ we consider the equivalency relation $x\sim y$ iff there is a scalar $a$ of magnitude one, $|a|=1$, so that
$y=ax$. Let $\hat{H}=H/\sim$ denote the set of equivalence classes. Note $\hat{H}\setminus\{0\}$ is equivalent to the cross-product between a real or complex projective
space ${\cal P}^{n-1}$ of dimension $n-1$ and the positive semiaxis $\R^{+}$.

In this paper we use $\hat{x}$ to denote the equivalency class of $x$ in $\hat{H}$. Nevertheless, for simplicity, $x$ is used in place of $\hat{x}$ when there is no ambiguity.

Let $\alpha$ and $\beta$ denote the nonlinear maps
\begin{equation}\label{eq:1}
\alpha: \hch\rightarrow\R^m~~,~~\alpha(x)=\left(|\ip{x}{f_k}|\right)_{1\leq k\leq m},
\end{equation}
\begin{equation}\label{eq:2}
\beta: \hch\rightarrow\R^m~~,~~\beta(x)=\left(|\ip{x}{f_k}|^2\right)_{1\leq k\leq m}.
\end{equation}

The {\em phase retrieval problem}, or the {\em phaseless reconstruction problem},
refers to analyzing when $\alpha$ (or $\beta$) is an injective map, and in this case to finding "good" left inverses.

The frame $\fc$ is said to be {\em phase retrievable} if the nonlinear map $\alpha$ (or $\beta$) is injective. 
In this paper we assume $\alpha$ and $\beta$ are injective maps (hence $\fc$ is phase retrievable). 
The problem is to analyze Lipschitz properties of these nonlinear maps, and then 
to extend the unique left inverse from the image of $\hat{H}$ through the nonlinear maps $\alpha$, $\beta$,
to the entire space $\R^m$ so that they remains Lipschitz continuous. 

A continuous map $f:(X,d_X)\rightarrow (Y,d_Y)$, defined between metric spaces $X$ and $Y$
with distances $d_X$ and $d_Y$ respectively, is Lipschitz continuous with Lipschitz constant $\Lip(f)$ if
\[ \Lip(f) := \sup_{x_1,x_2\in X} \frac{d_Y(f(x_1),f(x_2))}{d_X(x_1,x_2)} < \infty. \]
The map $f$ is called {\em bi-Lipschitz} with lower Lipschitz constant $a$ and upper Lipschitz constant $b$ if
for every $x_1,x_2\in X$,
\[ a\, d_X(x_1,x_2) \leq d_Y(f(x_1),f(x_2)) \leq b\, d_X(x_1,x_2). \]
Obviously the smallest upper Lischitz constant is $b=Lip(f)$. If $f$ is bi-Lipschitz then $f$ is injective. 

The space $\hch$ admits two classes of inequivalent metrics. We introduce and study them in detail in section \ref{sec2}. In particular consider the following two distances:
\begin{eqnarray}
D_2(x,y)  & = & \min_{\varphi}\norm{x-e^{i\varphi}y}_2 = \sqrt{\norm{x}^2+\norm{y}^2-2|\ip{x}{y}|} \\
d_1(x,y) & = & \norm{xx^*-yy^*}_1 = \sqrt{(\norm{x}^2+\norm{y}^2)^2-4|\ip{x}{y}|^2}
\end{eqnarray}
When the frame is phase retrievable the nonlinear maps $\alpha:(\hch,D_2)\rightarrow(\R^m,\norm{\cdot}_2)$
 and $\beta:(\hch,d_1)\rightarrow(\R^m,\norm{\cdot}_2)$ are shown to be bi-Lipschitz.  
This statement was
previously know for the map $\beta$ in the real and complex case (see \cite{Bal12a,Bal13a,BW13}), 
and for the map $\alpha$ in the real case only (see \cite{EM12,BW13,BCMN13}). 
In this paper we prove this statement for $\alpha$ in the complex case. 
Denote by $a_\alpha$ and $a_\beta$ the lower Lipschitz constants of $\alpha$ and $\beta$ respectively.
In this paper we prove also that there exist two Lipschitz continuous maps $\omega:(\R^m,\norm{\cdot}_2)\rightarrow
(\hch,D_2)$ and $\psi:(\R^m,\norm{\cdot}_2)\rightarrow(\hch,d_1)$ so that $\omega(\alpha(x))=x$ and
$\psi(\beta(x))=x$ for every $x\in \hc$. Furthermore the upper Lipschitz constants of these maps obey
$\Lip(\omega)\leq \frac{8.25}{a_\alpha}$ and $\Lip(\psi)\leq \frac{8.25}{a_\beta}$.
Surprisingly this shows the Lipschitz constant
of these left inverses are just a small factor larger than the minimal Lipschitz constants.
 Furthermore this factor is independent of dimension $n$ or number of frame vectors $m$.

The organization of the paper is as follows. Section \ref{sec2} introduces notations and presents the results for bi-Lipschitz properties. Section \ref{sec3} presents the results for the extension of the left inverse. Section \ref{sec4} contains the proof of these results.

\section{Notations and Bi-Lipschitz Properties}\label{sec2}

On the space $\hat{H}$ we consider two classes of metrics (distances) induced by corresponding distances on $H$ and $S^{1,0}(H)$ respectively:

\noindent 1. The class of {\em natural metrics}. For every $1\leq p \leq \infty$ and $x,y\in H$ define
\begin{equation}
D_p(\hat{x},\hat{y}) = \min_{|a|=1} \norm{x-ay}_p
\end{equation}
When no subscript is used, $\norm{\cdot}$ denotes the Euclidian norm, $\norm{\cdot}=\norm{\cdot}_2$.

\noindent 2. The class of {\em matrix norm induced metrics}. For every $1\leq p\leq \infty$ and $x,y\in H$ define
\begin{equation}\label{eq:dp}
d_p(\hat{x},\hat{y}) = \norm{xx^*-yy^*}_p = \left\{ \begin{array}{rcl}
\mbox{$\left(\sum_{k=1}^n (\sigma_k)^p \right) ^{1/p}$} & for & \mbox{$1\leq p\leq \infty$} \\
\mbox{$\max_{1\leq k\leq n} \sigma_k $} & for & \mbox{$p=\infty$}
\end{array} \right.
\end{equation}
where $(\sigma_k)_{1\leq k\leq n}$ are the singular values of the operator $xx^*-yy^*$, which is of rank at most 2.

Our choice in (\ref{eq:dp})
corresponds to the class of Schatten norms that extend to ideals of compact operators. In particular $p=\infty$ corresponds to the operator norm $\norm{\cdot}_{op}$ in
$\Sym(H)=\{T: H\rightarrow H~,~T=T^*\}$; $p=2$ corresponds to the Frobenius norm $\norm{\cdot}_{Fr}$ in $\Sym(H)$; $p=1$ corresponds to the nuclear norm $\norm{\cdot}_{*}$ in $\Sym(H)$:
\[ d_{\infty}(x,y) = \norm{xx^*-yy^*}_{op}~,~d_2(x,y)=\norm{xx^*-yy^*}_{Fr}~,~d_1(x,y)= \norm{xx^*-yy^*}_* \]
Note the Frobenius norm $\norm{T}_{Fr} = \sqrt{\trace(TT^*)}$ induces an Euclidian metric on $\Sym(H)$. In \cite{Bal13a} Lemma 3.7 we computed explicitly the eigenvalues of
 $\outp{x}{y}$. Based on these values, we can easily derive explicit expressions for these distances:
\[ d_{\infty}(\hat{x},\hat{y}) = \frac{1}{2}| \norm{x}^2-\norm{y}^2 | + \frac{1}{2}\sqrt{(\norm{x}^2+\norm{y}^2)^2-4|\ip{x}{y}|^2} \]
\[ d_{2}(x,y) = \sqrt{\norm{x}^4 + \norm{y}^4 - 2|\ip{x}{y}|^2 }\]
\[ d_1(x,y) = \sqrt{(\norm{x}^2+\norm{y}^2)^2-4|\ip{x}{y}|^2} \]
Since $xx^*-yy^*=\outp{u}{v}$ for $x=\frac{1}{2}(u+v)$ and $y=\frac{1}{2}(u-v)$.

To analyze the bi-Lipschitz properties, we define the following three types of Lipschitz bounds for $\alpha$. Note that the Lipschitz constants are square-roots of those constants.

\begin{enumerate}
\item The \emph{global lower} and \emph{upper Lipschitz bounds}, respectively:
\begin{equation*}
A_0 = \inf_{x,y \in \hat{H}} \frac{\norm{\alpha(x)-\alpha(y)}_2^2}{D_2(x,y)^2}~,
\qquad \qquad
B_0 = \sup_{x,y \in \hat{H}} \frac{\norm{\alpha(x)-\alpha(y)}_2^2}{D_2(x,y)^2}~;
\end{equation*}

\item The \emph{type I local lower} and \emph{upper Lipschitz bounds} at $z \in \hat{H}$, respectively:
\begin{equation*}
A(z) = \lim_{r \rightarrow 0} \inf_{\substack{x,y \in \hat{H}\\ D_2(x,z) < r\\ D_2(y,z) < r}} \frac{\norm{\alpha(x)-\alpha(y)}_2^2}{D_2(x,y)^2}~,
\qquad \qquad
B(z) = \sup_{r \rightarrow 0} \inf_{\substack{x,y \in \hat{H}\\ D_2(x,z) < r\\ D_2(y,z) < r}} \frac{\norm{\alpha(x)-\alpha(y)}_2^2}{D_2(x,y)^2}~;
\end{equation*}

\item The \emph{type II local lower} and \emph{upper Lipschitz bounds} at $z \in \hat{H}$, respectively:
\begin{equation*}
\tilde{A}(z) = \lim_{r \rightarrow 0} \inf_{\substack{x \in \hat{H}\\ D_2(x,z) < r}} \frac{\norm{\alpha(x)-\alpha(z)}_2^2}{D_2(x,z)^2}~,
\qquad \qquad
\tilde{B}(z) = \sup_{r \rightarrow 0} \inf_{\substack{x \in \hat{H}\\ D_2(x,z) < r}} \frac{\norm{\alpha(x)-\alpha(z)}_2^2}{D_2(x,y)^2}~.
\end{equation*}

\end{enumerate}

Similarly, we define the three types of Lipschitz constants for $\beta$.

\begin{enumerate}
\item The \emph{global lower} and \emph{upper Lipschitz bounds}, respectively:
\begin{equation*}
a_0 = \inf_{x,y \in \hat{H}} \frac{\norm{\beta(x)-\beta(y)}_2^2}{d_1(x,y)^2}~,
\qquad \qquad
b_0 = \sup_{x,y \in \hat{H}} \frac{\norm{\beta(x)-\beta(y)}_2^2}{d_1(x,y)^2}~;
\end{equation*}

\item The \emph{type I local lower} and \emph{upper Lipschitz bounds} at $z \in \hat{H}$, respectively:
\begin{equation*}
a(z) = \lim_{r \rightarrow 0} \inf_{\substack{x,y \in \hat{H}\\ d_1(x,z) < r\\ d_1(y,z) < r}} \frac{\norm{\beta(x)-\beta(y)}_2^2}{d_1(x,y)^2}~,
\qquad \qquad
b(z) = \lim_{r \rightarrow 0} \sup_{\substack{x,y \in \hat{H}\\ d_1(x,z) < r\\ d_1(y,z) < r}} \frac{\norm{\beta(x)-\beta(y)}_2^2}{d_1(x,y)^2}~;
\end{equation*}

\item The \emph{type II local lower} and \emph{upper Lipschitz bounds} at $z \in \hat{H}$, respectively:
\begin{equation*}
\tilde{a}(z) = \lim_{r \rightarrow 0} \inf_{\substack{y \in \hat{H}\\ d_1(x,z) < r}} \frac{\norm{\beta(x)-\beta(z)}_2^2}{d_1(x,z)^2}~,
\qquad \qquad
\tilde{b}(z) = \lim_{r \rightarrow 0} \sup_{\substack{y \in \hat{H}\\ d_1(x,z) < r}} \frac{\norm{\beta(x)-\beta(z)}_2^2}{d_1(x,z)^2}~.
\end{equation*}

\end{enumerate}

Note that due to homogeneity we have $A_0 = A(0)$, $B_0 = B(0)$, $a_0 = a(0)$, $b_0 = b(0)$. Also, for $z \neq 0$, we have $A(z) = A(z/\norm{z})$, $B(z) = B(z/\norm{z})$, $a(z) = a(z/\norm{z})$, $b(z) = b(z/\norm{z})$.

We analyze the bi-Lipschitz properties of $\alpha$ and $\beta$ by studying these constants.

\subsection{Bi-Lipschitz Properties of $\alpha$}\label{sec2_1}
The real case $H = \RR^n$ is studied in \cite{BW13}. We summarize the results as a theorem.

Recall that $\fc = \{f_1,\cdots,f_m\}$ is a frame in $H$ if there exist positive constants $A$ and $B$ for which
\begin{equation}\label{def:frame}
A \norm{x}^2 \leq \sum_{k=1}^m \absip{x}{f_k}^2 \leq B \norm{x}^2.
\end{equation}
We say $A$ [resp., $B$] is the optimal lower [resp., upper] frame bound if $A$ [resp., $B$] is the largest [resp., smallest] positive number for which the inequality (\ref{def:frame}) is satisfied.

For any index set $I \subset \{1,2,\cdots,m\}$, let $\cF[I] = \{ f_k, k \in I \}$ denote the frame subset indexed by $I$. Also, let $\sigma_1^2[I]$ and $\sigma_n^2[I]$ denote the upper and lower frame bound of set $\cF[I]$, respectively. That is,
\begin{equation*}
\sigma_1^2[I] = \l_{\max} \left( \sum_{k \in I} f_k f_k^{*} \right) \qquad \mbox{and} \qquad \sigma_n^2[I] = \l_{\min} \left( \sum_{k \in I} f_k f_k^{*} \right)~.
\end{equation*}

\begin{theorem}[\cite{BW13}]
Let $\cF \subset \RR^n$ be a phase retrievable frame for $\RR^n$. Let $A$ and $B$ denote its optimal lower and upper frame bound, respetively. Then


\begin{enumerate}
\item For every $0 \neq x \in \RR^n$, $A(x) = \sigma_n^2(\mbox{supp}(\alpha(x))$ where $\mbox{supp}(\alpha(x)) = \{ k, \ip{x}{f_k} \neq 0 \}$;

\item For every $x \in \RR^n$, $\tilde{A} = A$;

\item $A_0 = A(0) = \min_{I \subset \{1,2,\cdots,m \}} (\sigma_n^2[I]+\sigma_n^2[I^c]) > 0$;

\item For every $x \in \RR^n$, $B(x) = \tilde{B}(x) = B$;

\item $B_0 = B(0) = \tilde{B}(0) = B$.







\end{enumerate}
\end{theorem}

Now we consider the complex case $H = \cC^n$. We analyze the complex case by doing a realification first. Consider the $\RR$-linear map $\j: \cC^n \rightarrow \RR^{2n}$ defined by
\begin{equation*}
\j(z) = \left[ \begin{array}{c}
\mbox{real}(z)\\
\mbox{imag}(z)
\end{array} \right].
\end{equation*}
This realification is studied in detail in \cite{Bal13a}. We call $\j(z)$ the realification of $z$. For simplicity, in this paper we will denote $\xi = \j(x)$, $\eta = \j(y)$, $\zeta = \j(z)$, $\varphi = \j(f)$, $\delta = \j(d)$, respectively.

For a frame set $\cF = \{f_1,f_2,\cdots,f_m\}$, define the symmetric operator
\begin{equation*}
\Phi_k = \phik \varphi_k^T + J \phik \varphi_k^T J^T, \quad k = 1,2,\cdots,m.
\end{equation*}
where
\begin{equation}\label{def:J}
J = \left[ \begin{array}{cc}
0 & -I\\
I & 0
\end{array} \right]
\end{equation}
is a matrix in $\RR^{2n \times 2n}$.

Also, define $\SS: \RR^{2n} \rightarrow \Sym(\RR^{2n})$ by
\begin{equation*}
\SS(\xi) = \sum_{k: \Phi_k \xi \neq 0} \frac{1}{\ip{\Phi_k \xi}{\xi}} \Phi_k \xi \xi^T \Phi_k.
\end{equation*}
We have the following result:

\begin{theorem}\label{thm:alpha}
Let $\cF \subset \cC^n$ be a phase retrievable frame for $\cC^n$. Let $A$ and $B$ denote its optimal lower and upper frame bound, respetively. For any $z \in \cC^n$, let $\zeta = \j(z)$ be its realification. Then

\begin{enumerate}
\item For every $0 \neq z \in \cC^n$, $A(z) = \l_{2n-1}(\SS(\zeta))$ \label{item:2-2-1};

\item $A_0 = A(0) > 0$ \label{item:2-2-2};

\item For every $z \in \cC^n$, $\tilde{A}(z) = \l_{2n-1} \left( \SS(\zeta)+\sum_{k:\ip{z}{f_k}=0} \Phi_k \right)$ \label{item:2-2-3};

\item $\tilde{A}(0) = A$ \label{item:2-2-4};

\item For every $z \in \cC^n$, $B(z) = \tilde{B}(z) = \l_1 \left( \SS(\zeta)+\sum_{k:\ip{z}{f_k}=0} \Phi_k \right)$ \label{item:2-2-5};

\item $B_0 = B(0) = \tilde{B}(0) = B$ \label{item:2-2-6}.

\end{enumerate}
\end{theorem}

\subsection{Bi-Lipschitz Properties of $\beta$}\label{sec2_2}

The nonlinear map $\beta$ naturally induces a linear map between the space $\Sym(H)$ of symmetric operators on $H$ and $\R^m$:
\begin{equation}
\label{eq:A}
\A:\Sym(H)\rightarrow\R^m~~,~~\A(T)=(\ip{Tf_k}{f_k})_{1\leq k\leq m}
\end{equation}
This linear map has first been observed in \cite{BBCE07} and it has been exploited successfully in various paprs e.g. \cite{Bal2010,CSV12,Bal12a}.

Let $S^{p,q}(H)$ denote the set of symmetric operators that have at most $p$ strictly positive eigenvalues and $q$ strictly negative eigenvalues.
In particular $S^{1,0}(H)$ denotes the set of non-negative symmetric operators of rank at most one:
\begin{equation} S^{1,0}(H) =\{\, xx^*,~~ x\in H\, \} \end{equation}
where $x^*:H\rightarrow\C$, $x^*(y)=\ip{y}{x}$ is the dual map associated to vector $x$.
In \cite{Bal13a} we studied in more depth geometric and analytic properties of this set. In particular note
$\beta(x)=\A(\outp{x}{x})$ where
\begin{equation} \outp{x}{y}=\frac{1}{2}(xy^*+yx^*)  \end{equation}
denotes the symmetric outer product between vectors $x$ and $y$.
The map $\beta$ is injective if and only if $\A$ restricted to $S^{1,0}(H)$ is injective.

In previous papers \cite{Bal13a,BW13} we showed the following necessary and sufficient conditions for a frame to give phase retrieval.

\begin{theorem}[\cite{Bal13a,BW13}]
\label{equiv}
The following are equivalent:
\begin{enumerate}
\item The frame $\fc$ is phase retrievable;
\item $ker(\ac)\cap S^{1,1}(H) = \{0\}$;
\item There is a constant $a_0>0$ so that for every $u,v\in H$
\begin{equation}
\label{eq:lowerbound}
\frac{1}{2}\sum_{k=1}^m \left| \ip{u}{f_k}|\ip{f_k}{v}+\ip{v}{f_k}\ip{f_k}{u} \right|^2 \geq a_0 \left[ \norm{u}^2\norm{v}^2-(\imag(\ip{u}{v}))^2 \right]
\end{equation}
\item There is a constant $a_0>0$ so that for every $x,y\in H$,
\begin{equation}
\label{eq:a0}
\norm{\beta(x)-\beta(y)} ^2 \geq a_0 (d_1 (x,y))^2
\end{equation}
\end{enumerate}
\end{theorem}

In \cite{Bal13a} we also showed a theorem that can be restated as follows:

\begin{theorem}
\label{upperlip}
If $\fc$ is phase retrievable, then there is a constant $b_0 > 0$ such that for every $x, y \in H$,
\begin{equation}
\label{eq:b0given}
\norm{\beta(x)-\beta(y)} ^2 \leq b_0 (d_1 (x,y))^2
\end{equation}
where $b_0$ is given by
\begin{equation}
b_0 = \max_{\norm{x}=\norm{y}=1} \sum_{k=1}^{m} \left(\real\left(\ip{x}{f_k} \ip{f_k}{y} \right)\right)^2 = \max_{\norm{x} = 1} \sum_{k=1}^m \absip{x}{f_k}^4 = \norm{T}_{B(l^2,l^4)}^4 ~.
\end{equation}
In the last expression $T: H \rightarrow \cC^m$ is the analysis operator defined by $x \mapsto \left( \ip{x}{f_k} \right)_{k=1}^m$.
\end{theorem}

\begin{remark}
An upper bound of $b_0$ is given by
\begin{equation}
b_0  \leq B\left(\max_{1 \leq k \leq m} \norm{f_k} \right)^2 \leq B^2
\end{equation}
where $B$ is the upper frame bound of $\fc$.
\end{remark}

We give an expression of the local Lipshitz bounds as well. Define $\cR: \RR^{2n} \rightarrow \Sym(\RR^{2n})$ by
\begin{equation}
\cR(\xi) = \sum_{k=1}^m \Phi_k \xi \xi^T \Phi_k~.
\end{equation}

\begin{theorem}
\label{thm:beta}
Let $\cF$ be a phase retrievable frame for $H = \cC^n$. For every $0 \neq z \in H$, let $\zeta = \j(z)$ denote the realification of $z$. Then
\begin{enumerate}
\item $a(z) = \tilde{a}(z) = \l_{2n-1}(\cR(\zeta))/\norm{\zeta}^2$;
\item $b(z) = \tilde{b}(z) = \l_{1}(\cR(\zeta))/\norm{\zeta}^2$;
\item \textup{(\cite{Bal13a})} $a(0) = a_0 = \min_{\norm{\zeta}=1} \l_{2n-1} \left( \cR(\zeta) \right)$;
\item $\tilde{a}(0) = \min_{\norm{x} = 1} \sum_{k=1}^m \absip{x}{f_k}^4$;
\item $b(0) = \tilde{b}(0) = b_0$;
\end{enumerate}

\end{theorem}

\section{Extension of the Inverse Map}\label{sec3}

All metrics $D_p$ and $d_p$ induce the same topology as shown in the following result.
\begin{proposition}
\label{lem1}
\begin{enumerate}
\item For each $1\leq p\leq \infty$, $D_p$ and $d_p$ are metrics (distances) on $\hat{H}$.

\item $(D_p)_{1\leq p\leq\infty}$ are equivalent metrics, that is each $D_p$ induces the same topology on $\hat{H}$ as $D_1$.
Additionally, for every $1\leq p,q\leq\infty$
the embedding $i:(\hat{H},D_p)\rightarrow (\hat{H},D_q)$, $i(x)=x$,  is Lipschitz with Lipschitz constant
\begin{equation}\label{eq:L1}
L^D_{p,q,n}=\max(1,n^{\frac{1}{q}-\frac{1}{p}}).
\end{equation}

\item For $1\leq p,q\leq\infty$, $(d_p)_{1\leq p\leq\infty}$ are equivalent metrics, that is each $d_p$ induces the same topology on $\hat{H}$ as $d_1$.
Additionally, for every $1\leq p,q\leq\infty$
the embedding $i:(\hat{H},d_p)\rightarrow (\hat{H},d_q)$, $i(x)=x$,  is Lipschitz with Lipschitz constant
\begin{equation}\label{eq:L2}
L^d_{p,q,n}=\max(1,2^{\frac{1}{q}-\frac{1}{p}}).
\end{equation}

\item The identity map $i:(\hat{H},D_p)\rightarrow (\hat{H},d_p)$, $i(x)=x$, is continuous with continuous inverse.
However it is not Lipschitz, nor is its inverse.

\item The metric space $(\hat{H},D_p)$ is Lipschitz isomorphic to $S^{1,0}(H)$ endowed with Schatten norm $\norm{\cdot}_p$. The isomorphism is given by the map
\begin{equation}\label{eq:kappa_a}
\kappa_\alpha : \hat{H} \rightarrow S^{1,0}(H)~~,~~\kappa_\alpha(x)=\left\{ 
\begin{array}{ccc} \mbox{$\frac{1}{\norm{x}}xx^*$} & if & \mbox{$x\neq 0$} \\
 0 & if & x = 0 \end{array} \right..
\end{equation}
The embedding $\kappa_\alpha$ is bi-Lipschitz with the lower Lipschitz constant $\min(2^{\frac{1}{2}-\frac{1}{p}},n^{\frac{1}{p}-\frac{1}{2}})$ and the upper Lipschitz constant $\sqrt{2}\max(n^{\frac{1}{2}-\frac{1}{p}},2^{\frac{1}{p}-\frac{1}{2}})$. In particular for $p=2$ the lower Lipschitz constant is 1 and
the upper Lipschitz constant is $\sqrt{2}$.

\item The metric space $(\hat{H},d_p)$ is isometrically isomorphic to $S^{1,0}(H)$ endowed with Schatten norm $\norm{\cdot}_p$.
The isomorphism is given by the map
\begin{equation}\label{eq:kappa_b}
\kappa_\beta: \hat{H} \rightarrow S^{1,0}(H)~~,~~\kappa_\beta(x)= xx^*.
\end{equation}
In particular the metric space $(\hat{H},d_1)$ is isometrically isomorphic to $S^{1,0}(H)$ endowed with the nuclear norm $\norm{\cdot}_1$.

\item The nonlinear map $\iota: (\hat{H},D_p) \rightarrow (\hat{H},d_p)$ defined by
\begin{equation}
\iota(x) = \left\{ 
\begin{array}{ccc} \mbox{$\frac{x}{\sqrt{\norm{x}}}$} & if & \mbox{$x\neq 0$} \\
 0 & if & x = 0 \end{array} \right.
\end{equation}
is bi-Lipschitz with the lower Lipschitz constant $\min(2^{\frac{1}{2}-\frac{1}{p}}, n^{\frac{1}{p}-\frac{1}{2}})$ and the upper Lipschitz constant $\sqrt{2}\max(n^{\frac{1}{2}-\frac{1}{p}},2^{\frac{1}{p}-\frac{1}{2}})$.

\end{enumerate}
\end{proposition}

\begin{remark}

\begin{enumerate}
\item Note the Lipschitz bound $L^D_{p,q,n}$ is equal to the operator norm of the identity between $(\C^n,\norm{\cdot}_p)$ and $(\C^n,\norm{\cdot}_q)$:
$L^D_{p,q,n}=\norm{I}_{l^p(\C^n)\rightarrow l^q(\C^n)}$.

\item Note the equality $L^d_{p,q,n}=L^D_{p,q,2}$.
\end{enumerate}

\end{remark}

The results in Section \ref{sec2}, together with the previous proposition, show that frame $\fc$ is phase retrievable then the nonlinear map (\ref{eq:1}) [resp., (\ref{eq:2})]
is bi-Lipschitz between metric spaces $(\hat{H},D_p)$ [resp., $(\hat{H},d_p)$] and $(\R^m,\norm{\cdot}_q)$. Recall that the Lipschitz constants between $(\hat{H},D_2)$ [resp. $(\hat{H},d_1)$] and
$(\R^m,\norm{\cdot}=\norm{\cdot}_2)$ are given by $\sqrt{A_0}$ [resp., $\sqrt{a_0}$] and $\sqrt{B_0}$ [resp., $\sqrt{b_0}$]:
\begin{equation}\label{eq:lipr1}
\sqrt{A_0} D_2(x,y) \leq \norm{\alpha(x)-\alpha(y)} \leq \sqrt{B_0} D_2(x,y)
\end{equation}
\begin{equation}\label{eq:lipr2}
\sqrt{a_0} d_1(x,y) \leq \norm{\beta(x)-\beta(y)} \leq \sqrt{b_0} d_1(x,y)
\end{equation}
Clearly the inverse map defined on the range of $\alpha$ [resp., $\beta$] from metric space $(\alpha(\hat{H}),\norm{\cdot})$ [resp., $(\beta(\hat{H}),\norm{\cdot})$] to $(\hat{H},D_2)$ [resp., $(\hat{H},d_1)$]:
\begin{equation}
\label{eq:inv1}
\tilde{\omega} : \alpha(\hat{H})\subset\R^m \rightarrow \hat{H}~~,~~\tilde{\omega}(c) = x~~{\rm if}~\alpha(x)=c
\end{equation}
\begin{equation}
\label{eq:inv2}
\tilde{\psi} : \beta(\hat{H})\subset\R^m \rightarrow \hat{H}~~,~~\tilde{\psi}(c) = x~~{\rm if}~\beta(x)=c
\end{equation}
 is Lipschitz with Lipschitz constant $\frac{1}{\sqrt{A_0}}$ [resp., $\frac{1}{\sqrt{a_0}}$].
In this paper we prove that both $\tilde{\omega}$ and $\tilde{\psi}$ can be extended to the entire $\R^m$ as a Lipschitz map with Lipschitz constant that increases by a small factor.

The precise statement is given in the following Theorem which is the main result of this paper.

\begin{theorem}\label{maintheo}
Let $\fc=\{f_1,\ldots,f_m\}$ be a phase retrievable frame for the $n$ dimensional Hilbert space $H$, and let $\alpha, \beta :\hat{H}\rightarrow\R^m$ denote the
injective nonlinear analysis map as defined in (\ref{eq:1}) and (\ref{eq:2}). Let $A_0$ and $a_0$ denote the positive constant as in (\ref{eq:lipr1}) and (\ref{eq:lipr2}).
Then 

\begin{enumerate}
\item
there exists a Lipschitz continuous function $\omega:\R^m\rightarrow\hat{H}$ so that $\omega(\alpha(x))=x$ for all $x\in\hat{H}$.
For any $1\leq p,q\leq \infty$,  $\omega$ has an upper Lipschitz constant $\Lip(\omega)_{p,q}$ between $(\R^m,\norm{\cdot}_p)$ and $(\hat{H},D_q)$ bounded by:
\begin{equation}\label{eq:Lpq1}
\Lip(\omega)_{p,q} \leq \left\{ \begin{array}{cc}
\mbox{$\frac{3\sqrt{2}+4}{\sqrt{A_0}} \cdot 2^{\frac{1}{q}-\frac{1}{2}}  \cdot \max(1,m^{\frac{1}{2}-\frac{1}{p}}) $} & \mbox{for $q\leq 2$}\\
\mbox{$\frac{3\sqrt{2}+2^{\frac{3}{2}+\frac{1}{q}}}{\sqrt{A_0}} \cdot n ^{\frac{1}{2}-\frac{1}{q}} \cdot \max(1,m^{\frac{1}{2}-\frac{1}{p}}) $} & \mbox{for $q> 2$}
\end{array} \right.
\end{equation}
Explicitly this means: for $q\leq 2$ and for all $c,d\in\R^m$:
\begin{equation}
\label{eq:Lippq01}
D_q(\omega(c),\omega(d)) \leq \frac{3\sqrt{2}+4}{\sqrt{A_0}} \cdot 2^{\frac{1}{q}-\frac{1}{2}}  \cdot \max(1,m^{\frac{1}{2}-\frac{1}{p}})\norm{c-d}_p
\end{equation}
whereas for $q>2$ and for all $c,d\in\R^m$:
\begin{equation}
\label{eq:Lippq02}
D_q(\omega(c),\omega(d)) \leq \frac{3\sqrt{2}+2^{\frac{3}{2}+\frac{1}{q}}}{\sqrt{A_0}} \cdot n ^{\frac{1}{2}-\frac{1}{q}} \cdot \max(1,m^{\frac{1}{2}-\frac{1}{p}}) \norm{c-d}_p
\end{equation}
In particular, for $p=2$ and $q=2$ its Lipschitz constant $\Lip(\omega)_{2,2}$
bounded by $\frac{4+3\sqrt{2}}{\sqrt{a_0}}$:
\begin{equation}\label{eq:inv21}
D_2(\omega(c),\omega(d)) \leq \frac{4+3\sqrt{2}}{\sqrt{a_0}} \norm{c-d}
\end{equation}

\item
there exists a Lipschitz continuous function $\psi:\R^m\rightarrow\hat{H}$ so that $\psi(\beta(x))=x$ for all $x\in\hat{H}$.
For any $1\leq p,q\leq \infty$,  $\psi$ has an upper Lipschitz constant $\Lip(\psi)_{p,q}$ between $(\R^m,\norm{\cdot}_p)$ and $(\hat{H},d_q)$ bounded by:
\begin{equation}\label{eq:Lpq2}
\Lip(\psi)_{p,q} \leq \left\{ \begin{array}{cc}
\mbox{$\frac{3+2\sqrt{2}}{\sqrt{a_0}} \cdot 2^{\frac{1}{q}-\frac{1}{2}}  \cdot \max(1,m^{\frac{1}{2}-\frac{1}{p}})$} & \mbox{for $q\leq 2$}\\
\mbox{$\frac{3+2^{1+\frac{1}{q}}}{\sqrt{a_0}} \max(1,m^{\frac{1}{2}-\frac{1}{p}})$} & \mbox{for $q> 2$}
\end{array} \right.
\end{equation}
Explicitly this means: for $q\leq 2$ and for all $c,d\in\R^m$:
\begin{equation}
\label{eq:Lippq}
d_q(\psi(c),\psi(d)) \leq \frac{3+2\sqrt{2}}{\sqrt{a_0}} \cdot 2^{\frac{1}{q}-\frac{1}{2}}  \cdot \max(1,m^{\frac{1}{2}-\frac{1}{p}}) \norm{c-d}_p
\end{equation}
whereas for $q>2$ and for all $c,d\in\R^m$:
\begin{equation}
\label{eq:Lippq2}
d_q(\psi(c),\psi(d)) \leq \frac{3+2^{1+\frac{1}{q}}}{\sqrt{a_0}} \max(1,m^{\frac{1}{2}-\frac{1}{p}}) \norm{c-d}_p
\end{equation}
In particular, for $p=2$ and $q=1$ its Lipschitz constant $\Lip(\psi)_{2,1}$
bounded by $\frac{4+3\sqrt{2}}{\sqrt{a_0}}$:
\begin{equation}\label{eq:inv22}
d_1(\psi(c),\psi(d)) \leq \frac{4+3\sqrt{2}}{\sqrt{a_0}}\norm{c-d}
\end{equation}

\end{enumerate}
\end{theorem}

The proof of Theorem \ref{maintheo}, presented in Section \ref{sec3}, requires construction of a special Lipschitz map. We believe this particular result is interesting in itself and
may be used in other constructions. This construction is given in \cite{BZ15a} for the case $p=2$. Here we consider a general $p$ and give a better bound for the Lipschitz constant. We state it as a lemma.
\begin{lemma}
\label{prop:2}
Consider the spectral decomposition of any self-adjoint operator $A$ in $\Sym(H)$, say $A=\sum_{k=1}^d \lambda_{m(k)}P_k$, where
$\lambda_1\geq \lambda_2\geq\cdots\geq\lambda_n$ are the $n$ eigenvalues including multiplicities, and $P_1$,...,$P_d$ are the orthogonal projections
associated to the $d$ distinct eigenvalues. Additionally, $m(1)=1$ and $m(k+1)=m(k)+r(k)$, where $r(k)=rank(P_k)$ is the multiplicity of eigenvalue $\lambda_{m(k)}$.
Then the map
\begin{equation}
\label{eq:pi}
\pi: \Sym(H) \rightarrow S^{1,0}(H)~~,~~\pi(A)=(\lambda_1-\lambda_2)P_1
\end{equation}
satisfies the following two properties: 

\begin{enumerate}

\item for $1\leq p\leq \infty$, it is Lipschitz continuous from $(\Sym(H),\norm{\cdot}_{p})$ to $(S^{1,0}(H),\norm{\cdot}_{p})$
with Lipschitz constant less than or equal to $3+2^{1+\frac{1}{p}}$;  

\item $\pi(A)=A$ for all $A\in S^{1,0}(H)$.
\end{enumerate}
\end{lemma}

\begin{remark}
Numerical experiments suggest the Lipschitz constant of $\pi$ is smaller than 5 for $p=\infty$.
On the other hand it cannot be smaller than 2 as the following example shows.
\end{remark}

\begin{example}
If $A = \begin{pmatrix}
1 & 0 \\
0 & 1 \\\end{pmatrix}$,
$B = \begin{pmatrix}
2 & 0 \\
0 & 0 \\\end{pmatrix}$, then $\pi(A)=\begin{pmatrix}
0 & 0 \\
0 & 0 \\\end{pmatrix}$ and
$\pi(B) = \begin{pmatrix}
2 & 0 \\
0 & 0 \\\end{pmatrix}$.
Here we have $\norm{\pi(A)-\pi(B)}_{\infty}=2$ and $\norm{A-B}_{\infty}=1$. Thus for this example $\norm{\pi(A)-\pi(B)}_{\infty}=2\norm{A-B}_{\infty}$.
\end{example}

It is unlikely to obtain an isometric extension in Theorem \ref{maintheo}. Kirszbraun theorem \cite{WelWil75} gives a sufficient condition for
isometric extensions of Lipschitz maps.
The theorem states that isometric extensions are possible when the pair of metric spaces satisfy the Kirszbraun property, or the K property:
\begin{definition}
The Kirszbraun Property (K): Let $X$ and $Y$ be two metric spaces with metric $d_x$ and $d_y$ respectively. $(X,Y)$ is said to have Property (K) if for any pair of families of closed balls $\{B(x_i,r_i): i \in I\}$, $\{B(y_i,r_i): i \in I\}$, such that $d_y(y_i,y_j) \leq d_x(x_i,x_j)$ for each $i,j \in I$, it holds that $\bigcap B(x_i,r_i) \neq \emptyset \Rightarrow \bigcap B(y_i,r_i) \neq \emptyset$.
\end{definition}

If $(X,Y)$ has Property (K), then by Kirszbraun's Theorem we can extend a Lipschitz mapping defined on a subspace of $X$ to a Lipschitz mapping defined on $X$ while maintaining the Lipschitz constant. Unfortunately, if we consider $(X,d_X)=(\R^m,\norm{\cdot})$ and $Y=\hat{H}$, Property (K) does not hold for either $D_p$ or $d_p$.

\emph{Property (K) does not hold for $\hat{H}$ with norm $D_p$.} Specifically, $(\mathbb{R}^m,\mathbb{R}^n/\sim)$ does not have Property K.
\begin{example}
We give a counterexample for $m=n=2, p=2$: Let $\tilde{y}_1 = (3,1)$, $\tilde{y}_2 = (-1,1)$, $\tilde{y}_3 = (0,1)$ be the representatives of
three points $y_1$, $y_2$, $y_3$ in $\mathbb{R}^2/\sim$. Then $D_2(y_1, y_2) = 2\sqrt{2}$, $D_2(y_2, y_3) = 1$ and
$D_2(y_1, y_3) = 3$. Consider $x_1 = (0,0)$, $x_2 = (0,-2\sqrt{2})$, $x_3 = (-1,-2\sqrt{2})$ in $\mathbb{R}^2$ with the Euclidean distance,
then we have $\norm{x_1- x_2} = 2\sqrt{2}$, $\norm{x_2- x_3} = 1$ and $\norm{x_1- x_3} = 3$. For $r_1=\sqrt{6}$, $r_2=2-\sqrt{2}$, $r_3=\sqrt{6}-\sqrt{3}$,
we see that $(1-\sqrt{2},1+\sqrt{2}) \in \bigcap_{i=1}^{3}B(x_i,r_i)$ but $\bigcap_{i=1}^{3}B(y_i,r_i) = \emptyset$. To see $\bigcap_{i=1}^{3}B(y_i,r_i) = \emptyset$,
it suffices to look at the upper half plane in $\mathbb{R}^2$. If we look at the upper half plane $H$, then $B(y_1,r_1)$ becomes the union of two parts,
namely $B(\tilde{y}_1,r_1) \cup H$ and $B(-\tilde{y}_1,r_1) \cup H$, and $B(y_i,r_i)$ becomes $B(\tilde{y}_i,r_i)$ for $i=2$, $3$.
But $(B(\tilde{y}_1,r_1) \cup H) \cap B(\tilde{y}_2,r_2) = \emptyset$ and $(B(-\tilde{y}_1,r_1) \cup H) \cap B(\tilde{y}_3,r_3) = \emptyset$.
So we obtain that $\bigcap_{i=1}^{3}B(y_i,r_i) = \emptyset$.
\end{example}

\emph{Property (K) does not hold for $\hat{H}$ with norm $d_p$.} Specifically, $(\R^m,\mathbb{C}^n/\sim)$ does not have Property K. The following example is given in \cite{BZ15a}.
\begin{example}[\cite{BZ15a}]
Let $m$ be any positive integer and $n = 2$, $p = 2$. We want to show that $(X,Y) = (\RR^m,\cC^n / \sim)$ does not have Property (K). Let $\tilde{y}_1 = (1,0)$ and $\tilde{y}_2 = (0,\sqrt{3})$ be representitives of $y_1$, $y_2 \in Y$, respectively. Then $d_1(y_1,y_2) = 4$. Pick any two points $x_1$, $x_2$ in $X$ with $\norm{x_1-x_2} = 4$. Then $B(x_1,2)$ and $B(x_2,2)$ intersect at $x_3 = (x_1+x_2)/2 \in X$. It suffices to show that the closed balls $B(y_1,2)$ and $B(y_2,2)$ have no intersection in $H$. Assume on the contrary that the two balls intersect at $y_3$, then pick a representive of $y_3$, say $\tilde{y}_3 = (a,b)$ where $a$, $b \in \mathbb{C}$. It can be computed that
\begin{equation}
\label{ieq1}
d_1(y_1,y_3) =  \abs{a}^4+\abs{b}^4-2\abs{a}^2+2\abs{b}^2+2\abs{a}^2\abs{b}^2+1
\end{equation}
and
\begin{equation}
\label{ieq2}
d_1(y_2,y_3) =  \abs{a}^4+\abs{b}^4+6\abs{a}^2-6\abs{b}^2+2\abs{a}^2\abs{b}^2+9
\end{equation}
Set $d_1(y_1,y_3) = d_1(y_2,y_3) = 2$. Take the difference of the right hand side of (\ref{ieq1}) and (\ref{ieq2}), we have $\abs{b}^2-\abs{a}^2 = 1$ and thus $\abs{b}^2 \geq 1$. However, the right hand side of (\ref{ieq1}) can be rewritten as $(\abs{a}^2+\abs{b}^2-1)^2+4\abs{b}^2$, so $d_1(y_1,y_3) = 2$ would imply that $\abs{b}^2 \leq 1/2$. This is a contradiction.
\end{example}

\begin{remark}
Using nonlinear functional analysis language (\cite{BenLin00}) Lemma \ref{prop:2} can be restated by saying that $S^{1,0}(H)$ is a 5-Lipschitz retract in $\Sym(H)$.
\end{remark}

\begin{remark}
The Lipschitz inversion results of Theorem \ref{maintheo} can easily be extended to systems of quadratic equations, not necessarily of rank-1 matrices from the phase retrieval model considered in this paper.
\end{remark}

\section{Proofs of results}\label{sec4}

\subsection{Proof of results in Section \ref{sec2}}\label{sec4_1}
We start by proving Theorem \ref{thm:alpha}.

\noindent{\bf Proof of Theorem \ref{thm:alpha}}

\begin{enumerate}

\item First we prove the following lemma.
\begin{lemma}\label{lem:local}
Fix $x \in \cC^n$ and $z \in \cC^n$. Let $\xi = \j(x)$ and $\zeta = \j(z)$ be their realifications, respectively. Let $\xi_0 \in \hat{\xi} := \{\j(\tilde{x}) \in \RR^{2n}: \tilde{x} \in \hat{x}\}$ be a point in the equivalency class that satisfies $D_2(x,z) = \norm{\xi_0-\zeta}$. Then it is necessary that 
\begin{equation}\label{eq:lemlocal1}
\ip{\xi_0}{J\zeta} = 0
\end{equation}
and
\begin{equation}\label{eq:lemlocal2}
\ip{\xi_0}{\zeta} \geq 0
\end{equation}
where $J$ is defined as in (\ref{def:J}).
\end{lemma}

\emph{Proof}: For $\theta \in [0,2\pi)$ define
\begin{equation}
U(\theta) := \cos(\theta)I + \sin(\theta)J ~.
\end{equation}
Then it is easy to compute that
\begin{equation}
\j(e^{i\theta}x) = U(\theta)\xi ~.
\end{equation}
Therefore,
\begin{equation*}
D_2(x,z) = \min_{\theta \in [0,2\pi)} \norm{U(\theta)\xi-\zeta}^2 = \norm{\xi}^2+\norm{\zeta}^2-2\max_{\theta \in [0,2\pi)} \ip{U(\theta)\xi}{\zeta} ~.
\end{equation*}
If $\ip{U(\theta)\xi}{\zeta}$ is constantly zero, then we are done. Otherwise, note that
\begin{equation}
\max_{\theta \in [0,2\pi)} \ip{U(\theta)\xi}{\zeta} = \left(\ip{\xi}{\zeta}^2+\ip{J\xi}{\zeta}^2 \right)^{\frac{1}{2}}
\end{equation}
and the maximum is achieved at $\theta = \theta_0$ if and only if
\begin{equation}
\cos(\theta_0) = \frac{\ip{\xi}{\zeta}}{\left(\ip{\xi}{\zeta}^2+\ip{J\xi}{\zeta}^2 \right)^{\frac{1}{2}}}
\end{equation}
and
\begin{equation}
\sin(\theta_0) = \frac{\ip{J\xi}{\zeta}}{\left(\ip{\xi}{\zeta}^2+\ip{J\xi}{\zeta}^2 \right)^{\frac{1}{2}}} ~.
\end{equation}
Now we can compute
\begin{equation*}
\begin{aligned}
\ip{\xi_0}{J\zeta} & = \ip{U(\theta_0)\xi}{J\zeta} \\
& = \cos(\theta_0) \ip{\xi}{J\zeta} + \sin(\theta_0) \ip{J\xi}{J\zeta} \\
& = \frac{\ip{\xi}{\zeta}}{\left(\ip{\xi}{\zeta}^2+\ip{J\xi}{\zeta}^2 \right)^{\frac{1}{2}}} \ip{\xi}{J\zeta} + \frac{\ip{J\xi}{\zeta}}{\left(\ip{\xi}{\zeta}^2+\ip{J\xi}{\zeta}^2 \right)^{\frac{1}{2}}} \ip{J\xi}{J\zeta} \\
& = \frac{\ip{\xi}{\zeta}}{\left(\ip{\xi}{\zeta}^2+\ip{J\xi}{\zeta}^2 \right)^{\frac{1}{2}}} \ip{-J\xi}{\zeta} + \frac{\ip{J\xi}{\zeta}}{\left(\ip{\xi}{\zeta}^2+\ip{J\xi}{\zeta}^2 \right)^{\frac{1}{2}}} \ip{\xi}{\zeta} \\
& = 0 ~.
\end{aligned}
\end{equation*}
So we get (\ref{eq:lemlocal1}). (\ref{eq:lemlocal2}) is obvious.

\emph{Q.E.D.}

Now we come back to the proof of the theorem. Denote
\begin{equation}
p(x,y) :=  \frac{\norm{\alpha(x)-\alpha(y)}^2}{D_2(x,y)^2}, \qquad x,y \in \cC^n, ~\hat{x} \neq \hat{y}.
\end{equation}
We can represent this quotient in terms of $\xi$ and $\eta$. It is easy to compute that
\begin{equation}\label{eq:pP}
p(x,y) = P(\xi,\eta) := \frac{\sum_{k=1}^m \ip{\Phi_k \xi}{\xi} + \ip{\Phi_k \eta}{\eta} - 2\sqrt{\ip{\Phi_k \xi}{\xi} \ip{\Phi_k \eta}{\eta}}}{\norm{\xi}^2 + \norm{\eta}^2 - 2\sqrt{\ip{\xi}{\eta}^2+\ip{\xi}{J\eta}^2}} ~.
\end{equation}

Fix $r > 0$. Take $\xi$, $\eta \in \RR^{2n}$ that satisfy $D_2(x,z) = \norm{\xi-\zeta} < r$ and $D_2(y,z) = \norm{\eta-\zeta} < r$. Let $\mu = (\xi+\eta)/2$ and $\nu = (\xi-\eta)/2$. Then $\norm{\nu} < r$. Note that for $r$ small enough we have that $\norm{\mu} > \norm{\nu}$ and that $\Phi_k \zeta \neq 0 \Rightarrow \Phi_k \mu \neq 0$. Thus
\begin{equation*}
\begin{aligned}
P(\xi,\eta) & = \frac{\sum_{k=1}^m \ip{\Phi_k (\mu+\nu)}{\mu+\nu} + \ip{\Phi_k (\mu-\nu)}{\mu-\nu} - 2\sqrt{\ip{\Phi_k (\mu+\nu)}{\mu+\nu} \ip{\Phi_k (\mu-\nu)}{\mu-\nu}}}{\norm{\mu+\nu}^2 + \norm{\mu-\nu}^2 - 2\sqrt{\ip{\mu+\nu}{\mu-\nu}^2+\ip{\mu+\nu}{J(\mu-\nu)}^2}} \\
& = \frac{\sum_{k=1}^m \ip{\Phi_k \mu}{\mu} + \ip{\Phi_k \nu}{\nu} - \sqrt{(\ip{\Phi_k \mu}{\mu} + \ip{\Phi_k \nu}{\nu})^2 - 4\ip{\Phi_k \mu}{\nu}^2}}{\norm{\mu}^2 + \norm{\nu}^2 - \sqrt{\norm{\mu}^4+\norm{\nu}^4-2\norm{\mu}^2\norm{\nu}^2+4\ip{\mu}{J\nu}^2}} \\
& \geq \frac{\sum_{k: \Phi_k \zeta \neq 0} \ip{\Phi_k \mu}{\mu} + \ip{\Phi_k \nu}{\nu} - \sqrt{(\ip{\Phi_k \mu}{\mu} + \ip{\Phi_k \nu}{\nu})^2 - 4\ip{\Phi_k \mu}{\nu}^2}}{\norm{\mu}^2 + \norm{\nu}^2 - \sqrt{\norm{\mu}^4+\norm{\nu}^4-2\norm{\mu}^2\norm{\nu}^2}} \\
& = \frac{\sum_{k: \Phi_k \zeta \neq 0} \ip{\Phi_k \mu}{\mu} + \ip{\Phi_k \nu}{\nu} - \sqrt{(\ip{\Phi_k \mu}{\mu} + \ip{\Phi_k \nu}{\nu})^2 - 4\ip{\Phi_k \mu}{\nu}^2}}{2\norm{\nu}^2} \\
& = \frac{\sum_{k: \Phi_k \zeta \neq 0} \ip{\Phi_k \mu}{\mu} + \ip{\Phi_k \nu}{\nu} - \ip{\Phi_k \mu}{\mu} \sqrt{\left( 1+\frac{\ip{\Phi_k \nu}{\nu}}{\ip{\Phi_k \mu}{\mu}} \right)^2 - 4\frac{\ip{\Phi_k \mu}{\nu}^2}{\ip{\Phi_k \mu}{\mu}^2}}}{2\norm{\nu}^2} \\
& = \frac{\sum_{k: \Phi_k \zeta \neq 0} \ip{\Phi_k \mu}{\mu} + \ip{\Phi_k \nu}{\nu} - \ip{\Phi_k \mu}{\mu} \sqrt{1 + 2\frac{\ip{\Phi_k \nu}{\nu}}{\ip{\Phi_k \mu}{\mu}} +\frac{\ip{\Phi_k \nu}{\nu}^2}{\ip{\Phi_k \mu}{\mu}^2} - 4\frac{\ip{\Phi_k \mu}{\nu}^2}{\ip{\Phi_k \mu}{\mu}^2}}}{2\norm{\nu}^2} \\
& = \frac{\sum_{k: \Phi_k \zeta \neq 0} \ip{\Phi_k \mu}{\mu} + \ip{\Phi_k \nu}{\nu} - \ip{\Phi_k \mu}{\mu} \left( 1 + \frac{\ip{\Phi_k \nu}{\nu}}{\ip{\Phi_k \mu}{\mu}} - 2\frac{\ip{\Phi_k \mu}{\nu}^2}{\ip{\Phi_k \mu}{\mu}^2} \right) + O(\norm{\nu}^4) }{2\norm{\nu}^2} \\
& = \sum_{k: \Phi_k \zeta \neq 0} \frac{\ip{\Phi_k \mu}{\nu}^2}{\ip{\Phi_k \mu}{\mu} \norm{\nu}^2} + O(\norm{\nu}^2) \\
& = \frac{1}{\norm{\nu}^2} \ip{\SS(\mu)\nu}{\nu} + O(\norm{\nu}^2) ~.
\end{aligned}
\end{equation*}
Note that 
\begin{equation}\label{eq:iplnorm}
\absip{J\mu}{\nu} = \abs{\ip{J\mu}{\nu} - \ip{J\zeta}{\nu}} \leq \norm{J\mu-J\zeta}\norm{\nu} = \norm{\mu-\zeta}\norm{\nu}
\end{equation}
since $\ip{J\zeta}{\nu} = 0$ by Lemma \ref{lem:local}. Also, $\norm{\mu-\zeta} < r$. Therefore,
\begin{equation*}
\norm{P_{J\mu} \nu} = \frac{\absip{J\mu}{\nu}}{\norm{J\mu}} = \frac{\absip{J\mu}{\nu}}{\norm{\mu}} \leq \frac{r\norm{\nu}}{\norm{\mu}}
\end{equation*}
and thus
\begin{equation*}
\norm{P_{J\mu}^{\perp} \nu}^2 \geq \left( 1-\frac{r^2}{\norm{\mu}^2} \right) \norm{\nu}^2 ~.
\end{equation*}
As a consequence, we have
\begin{equation*}
\begin{aligned}
P(\xi,\eta) & = \frac{\ip{\SS(\mu)P_{J\mu}^{\perp}\nu}{P_{J\mu}^{\perp}\nu}}{\norm{\nu}^2} + O(\norm{\nu}^2) \\
& \geq \frac{\ip{\SS(\mu)P_{J\mu}^{\perp}\nu}{P_{J\mu}^{\perp}\nu}}{\norm{P_{J\mu}^{\perp}\nu}^2} \left( 1-\frac{r^2}{\norm{\mu}^2} \right) + O(r^2) \\
& \geq \left( 1-\frac{r^2}{\norm{\mu}^2} \right) \l_{2n-1} \left( \SS(\mu) \right) + O(r^2) ~.
\end{aligned}
\end{equation*}
Take $r \rightarrow 0$, by the continuity of eigenvalues with respect to matrix entries we have that
\begin{equation}\label{eq:Ageql}
A(z) \geq \l_{2n-1}(\SS(\zeta)) ~.
\end{equation}
On the other hand, take $E_{2n-1}$ to be the unit-norm eigenvector correspondent to $\l_{2n-1}(\SS(\zeta))$. For each $r>0$, take $\xi = \zeta+\frac{r}{2}E_{2n-1}$ and $\eta = \zeta-\frac{r}{2}E_{2n-1}$. Then
\begin{equation*}
p(x,y) = P(\xi,\eta) = \l_{2n-1}(\SS(\zeta)) ~.
\end{equation*}
Hence
\begin{equation}
A(z) \leq \l_{2n-1}(\SS(\zeta)) ~.
\end{equation}
Together with (\ref{eq:Ageql}) we have
\begin{equation}
A(z) = \l_{2n-1}(\SS(\zeta)) ~.
\end{equation}

\item Assume on the contrary that $A_0 = 0$, then for any $N \in \NN$, there exist $x_N$, $y_N \in H$ for which
\begin{equation}\label{eq:cpt_a}
p(x_N,y_N) = \frac{\norm{\alpha(x_N)-\alpha(y_N)}^2}{D_2(x_N,y_N)^2} \leq \frac{1}{N}.
\end{equation}
Without loss of generality we assume that $\norm{x_N} \geq \norm{y_N}$ for each $N$, for otherwise we can just swap the role of $x_N$ and $y_N$. Also due to homogeneity we assume $\norm{x_N} = 1$. By compactness of the closed ball $\cB_{1}(0) = \{x \in H: \norm{x} \leq 1 \}$ in $H = \cC^n$, there exist convergent subsequences of $\{x_N\}_{N \in \NN}$ and $\{y_N\}_{N \in \NN}$, which to avoid overuse of notations we still denote as $\{x_N\}_{N \in \NN} \rightarrow x_0 \in H$ and $\{y_N\}_{N \in \NN} \rightarrow y_0 \in H$.


Since $\norm{x_0} = 1$ we have from (i) that $A(x_0) > 0$. Note that $D_2(x_N,y_N) \leq \norm{x_N} + \norm{y_N} \leq 2$, so by (\ref{eq:cpt_a}) we have $\norm{\alpha(x_N)-\alpha(y_N)} \rightarrow 0$. That is, $\norm{\alpha(x_0)-\alpha(y_0)} = 0$. By injectivity we have $x_0 = y_0$ in $\hat{H}$. By Proposition \ref{thm:alpha}(i), 
\begin{equation*}
p(x_N,y_N) \geq A(x_0)-1/N > 1/N
\end{equation*}
for $N$ large enough. This is a contradiction with (\ref{eq:cpt_a}).

\item The case $z = 0$ is an easy computation. We now present the proof for $z \neq 0$. First we consider $p(x,z) = P(\xi,\zeta)$ as defined in (\ref{eq:pP}). Fix $r>0$. Take $\xi \in \RR^{2n}$ that satisfy $D_2(x,z)= \norm{\xi-\zeta} < r$. Let $d = x-z$ and $\delta = \j(d) = \xi-\zeta$. Note that
\begin{equation*}
P(\xi,\zeta) = \frac{\sum_{k=1}^m \ip{\Phi_k \xi}{\xi} + \ip{\Phi_k \zeta}{\zeta} - 2\sqrt{\ip{\Phi_k \xi}{\xi} \ip{\Phi_k \zeta}{\zeta}}}{\norm{\xi}^2 + \norm{\zeta}^2 - 2\sqrt{\ip{\xi}{\zeta}^2+\ip{\xi}{J\zeta}^2}} ~.
\end{equation*}
The numerator is equal to
\begin{equation*}
\begin{aligned}
& \sum_{k=1}^{m} \ip{\Phi_k \zeta}{\zeta} + 2\ip{\Phi_k \zeta}{\delta} + \ip{\Phi_k \delta}{\delta} + \ip{\Phi_k \zeta}{\zeta} - 2 \sqrt{(\ip{\Phi_k \zeta}{\zeta} + 2\ip{\Phi_k \zeta}{\delta} + \ip{\Phi_k \delta}{\delta}) \cdot \ip{\Phi_k \zeta}{\zeta}} \\
= & \sum_{k:\Phi_k \zeta \neq 0} 2\ip{\Phi_k \zeta}{\zeta} + 2\ip{\Phi_k \zeta}{\delta} + \ip{\Phi_k \delta}{\delta} + 2 \ip{\Phi_k \zeta}{\zeta} [ 1 + \frac{\ip{\Phi_k \zeta}{\zeta}\ip{\Phi_k \zeta}{\delta} + \frac{1}{2}\ip{\Phi_k \zeta}{\zeta}\ip{\Phi_k \delta}{\delta}}{\ip{\Phi_k \zeta}{\zeta}^2} - \\
& \quad \frac{1}{8} \cdot \frac{4 \ip{\Phi_k \zeta}{\zeta}^2 \ip{\Phi_k \zeta}{\delta}^2}{\ip{\Phi_k \zeta}{\zeta}^4} + O(\norm{\delta}^3) ] + \sum_{k:\Phi_k \zeta = 0} \ip{\Phi_k \delta}{\delta} \\
= & \sum_{k:\Phi_k \zeta \neq 0} \frac{\ip{\Phi_k \zeta}{\delta}^2}{\ip{\Phi_k \zeta}{\zeta}} + \sum_{k:\Phi_k \zeta = 0} \ip{\Phi_k \delta}{\delta} + O(\norm{\delta}^3) ~;
\end{aligned}
\end{equation*}
the denominator is equal to
\begin{equation*}
\begin{aligned}
& 2 \norm{\zeta}^2 + \norm{\delta}^2 + 2\ip{\zeta}{\delta} - 2 \norm{\zeta}^2 \left( 1 +  \frac{\norm{\zeta}^2 \ip{\zeta}{\delta} + \frac{1}{2} \ip{\zeta}{\delta} + \frac{1}{2} \ip{J \zeta}{\delta}^2}{\norm{\zeta}^4} - \frac{4\norm{\zeta}^4 \ip{\zeta}{\delta}^2}{8\norm{\zeta}^8} + O(\norm{\delta}^3) \right) \\
= & \norm{\delta}^2 + O(\norm{\delta}^3)
\end{aligned}
\end{equation*}
in which we used Lemma \ref{lem:local} to get $\ip{J\zeta}{\delta} = 0$.

Take $r \rightarrow 0$, we see that
\begin{equation}
\tilde{A}(z) \geq \l_{2n-1} \left( \SS(\zeta)+\sum_{k:\ip{z}{f_k}=0} \Phi_k \right) ~.
\end{equation}
Let $\tilde{E}_{2n-1}$ be the unit-norm eigenvector correspondent to\ $\l_{2n-1} \left( \SS(\zeta)+\sum_{k:\ip{z}{f_k}=0} \Phi_k \right)$. Note that $\ip{J\zeta}{\tilde{E}_{2n-1}} = 0$ since $\SS(\zeta)J\zeta = 0$ and $\Phi_k J\zeta = J\Phi_k \zeta = 0$ for each $k$ with $\ip{z}{f_k} = 0$.
Take $\xi = \zeta + \frac{r}{2} \tilde{E}_{2n-1}$ for each $r$, we again also have
\begin{equation}
\tilde{A}(z) \leq \l_{2n-1} \left( \SS(\zeta)+\sum_{k:\ip{z}{f_k}=0} \Phi_k \right) ~.
\end{equation}
Therefore
\begin{equation}
\tilde{A}(z) = \l_{2n-1} \left( \SS(\zeta)+\sum_{k:\ip{z}{f_k}=0} \Phi_k \right) ~.
\end{equation}

\item Take $z=0$ in (iii).

\item $\tilde{B}(z)$ can be computed in a similar way as in (iii) (in particular, the expansion for $P(\xi,\zeta)$ is exactly the same). We compute $B(z)$. $B(0)$ is computed in \cite{BCMN13}, Lemma 16. Now we consider $z \neq 0$. Use the same notations as in (\ref{eq:pP}). Fix $r>0$. Again, take $\xi$, $\eta \in \RR^{2n}$ that satisfy $D_2(x,z) = \norm{\xi-\zeta} < r$ and $D_2(y,z) = \norm{\eta-\zeta} < r$. Let $\mu = (\xi+\eta)/2$ and $\nu = (\xi-\eta)/2$. Also let $\delta_1 = \xi-\zeta$ and $\delta_2 = \eta-\zeta$.
Recall that 
\begin{equation*}
\begin{aligned}
P(\xi,\eta) & = \frac{\sum_{k=1}^m \ip{\Phi_k \xi}{\xi} + \ip{\Phi_k \eta}{\eta} - 2\sqrt{\ip{\Phi_k \xi}{\xi} \ip{\Phi_k \eta}{\eta}}}{\norm{\xi}^2 + \norm{\eta}^2 - 2\sqrt{\ip{\xi}{\eta}^2+\ip{\xi}{J\eta}^2}} \\
& = \sum_{k=1}^m \frac{ \ip{\Phi_k \xi}{\xi} + \ip{\Phi_k \eta}{\eta} - 2\sqrt{\ip{\Phi_k \xi}{\xi} \ip{\Phi_k \eta}{\eta}}}{\norm{\xi}^2 + \norm{\eta}^2 - 2\sqrt{\ip{\xi}{\eta}^2+\ip{\xi}{J\eta}^2}} ~.
\end{aligned}
\end{equation*}
Now we compute it as $\sum_{k=1}^m = \sum_{k:\Phi_k\zeta \neq 0} + \sum_{k:\Phi_k\zeta = 0}$. Again
\begin{equation}\label{eq:localB1}
\begin{aligned}
& \sum_{k: \Phi_k \zeta \neq 0} \frac{ \ip{\Phi_k \xi}{\xi} + \ip{\Phi_k \eta}{\eta} - 2\sqrt{\ip{\Phi_k \xi}{\xi} \ip{\Phi_k \eta}{\eta}}}{\norm{\xi}^2 + \norm{\eta}^2 - 2\sqrt{\ip{\xi}{\eta}^2+\ip{\xi}{J\eta}^2}} \\
= & \sum_{k: \Phi_k \zeta \neq 0} \frac{\ip{\Phi_k \mu}{\mu} + \ip{\Phi_k \nu}{\nu} - \sqrt{(\ip{\Phi_k \mu}{\mu} + \ip{\Phi_k \nu}{\nu})^2 - 4\ip{\Phi_k \mu}{\nu}^2}}{\norm{\mu}^2 + \norm{\nu}^2 - \sqrt{\norm{\mu}^4+\norm{\nu}^4-2\norm{\mu}^2\norm{\nu}^2+4\ip{\mu}{J\nu}^2}}
\end{aligned}
\end{equation}
The computation for its numerator is the same as in (i). We get that the numerator is equal to
\begin{equation*}
2\ip{\SS(\mu)\nu}{\nu} + O(\norm{\nu}^4)~.
\end{equation*}
Since $\mu \neq 0$, the denominator is equal to
\begin{equation}\label{eq:localB2}
\begin{aligned}
& \norm{\mu}^2+\norm{\nu}^2-\norm{\mu}^2\sqrt{1+\frac{\norm{\nu}^4}{\norm{\mu}^4}-\frac{2\norm{\nu}^2}{\norm{\mu}^2}+\frac{4\ip{\mu}{J\nu}^2}{\norm{\mu}^4}} \\
= & \norm{\mu}^2+\norm{\nu}^2-\norm{\mu}^2 \left( 1-\frac{\norm{\nu}^2}{\norm{\mu}^2}+\frac{2\ip{\mu}{J\nu}^2}{\norm{\mu}^4} \right) + O(\norm{\nu}^4) \\
= & 2\norm{\nu}^2 - \frac{2\ip{J\mu}{\nu}^2}{\norm{\mu}^2} + O(\norm{\nu}^4) \\
= & 2\norm{\nu}^2 + O(\norm{\nu}^4) \quad \mbox{by (\ref{eq:iplnorm})}.
\end{aligned}
\end{equation}
Also we can compute using the denominator as above [note that $\nu = (\delta_1-\delta_2)/2$] that
\begin{equation}\label{eq:localB3}
\begin{aligned}
& \sum_{k: \Phi_k \zeta = 0} \frac{ \ip{\Phi_k \xi}{\xi} + \ip{\Phi_k \eta}{\eta} - 2\sqrt{\ip{\Phi_k \xi}{\xi} \ip{\Phi_k \eta}{\eta}}}{\norm{\xi}^2 + \norm{\eta}^2 - 2\sqrt{\ip{\xi}{\eta}^2+\ip{\xi}{J\eta}^2}} \\
= & \sum_{k: \Phi_k \zeta = 0} \frac{ \left( \norm{\Phi_k^{1/2}\delta_1}-\norm{\Phi_k^{1/2}\delta_2} \right)^2 }{\norm{\delta_1-\delta_2}^2 + O(\norm{\nu}^4)}
\end{aligned}
\end{equation}
Now put together (\ref{eq:localB1}), (\ref{eq:localB2}) and (\ref{eq:localB3}), we get
\begin{equation}
P(\xi,\eta) = \frac{\ip{\SS(\mu)\nu}{\nu}+O(\norm{\nu}^4)}{\norm{\nu}^2+O(\norm{\nu}^4)} + \sum_{k:\Phi_k \zeta = 0} \frac{ \left( \norm{\Phi_k^{1/2}\delta_1}-\norm{\Phi_k^{1/2}\delta_2} \right)^2 }{\norm{\delta_1-\delta_2}^2 + O(\norm{\nu}^4)} ~.
\end{equation}
Note that
\begin{equation*}
\left( \norm{\Phi_k^{1/2} \delta_1} - \norm{\Phi_k^{1/2} \delta_2} \right)^2 \leq \ip{\Phi_k(\delta_1-\delta_2)}{\delta_1-\delta_2}
\end{equation*}
since it is equivalent to
\begin{equation}\label{iq:C-S}
\ip{\Phi_k \delta_1}{\delta_1} \ip{\Phi_k \delta_2}{\delta_2} \geq \left( \ip{\Phi_k \delta_1}{\delta_2} \right)^2
\end{equation}
which is the Cauchy-Schwarz inequality. Therefore we have that
\begin{equation}
P(\xi,\eta) \leq \frac{\ip{ ( \SS(\mu)+\sum_{k:\Phi_k \zeta = 0} \Phi_k ) \nu}{\nu} + O(\norm{\nu}^4)}{\norm{\nu}^2 + O(\norm{\nu}^4)} \leq \l_1 \left( \SS(\mu)+\sum_{k:\Phi_k \zeta = 0} \Phi_k \right) + O(r^2) ~.
\end{equation}
Take $r \rightarrow 0$ we have that
\begin{equation*}
B(z) \leq \l_1 \left( \SS(\zeta)+\sum_{k:\Phi_k \zeta = 0} \Phi_k \right)~.
\end{equation*}

Again we get the other direction of the above inequality by taking $\xi = \zeta + \frac{r}{2} E_1$ and $\eta = \zeta -  \frac{r}{2} E_1$ for each $r > 0$ where $E_1$ is the unit-norm eigenvector correspondent to $\l_1 \left( \SS(\zeta)+\sum_{k:\ip{z}{f_k}=0} \Phi_k \right)$. Note that for each $r$, the equality in (\ref{iq:C-S}) holds for this pair of $\xi$ and $\eta$.

\item Take $z=0$ in (v).

\end{enumerate}

Now we prove Theorem \ref{thm:beta}.

\noindent{\bf Proof of Theorem \ref{thm:beta}}

Only the first two parts are nontrivial. We prove them as follows.

Fix $z \in \cC^n$. Take $x = z+d_1$ and $y = z+d_2$ with $\norm{d_1}<r$ and $\norm{d_2}<r$ for $r$ small. Let $u = x+y = 2z+d_1+d_2$ and $v = x-y = d_1-d_2$. Let $\mu = 2\zeta+\delta_1+\delta_2 \in \RR^{2n}$ and $\nu = \delta_1-\delta_2 \in \RR^{2n}$ be the realification of $u$ and $v$, respectively. Define
\begin{equation}
\rho(x,y) = \frac{\norm{\beta(x)-\beta(y)}^2}{d_1(x,y)^2} ~.
\end{equation}
By the same computation as in \cite{Bal13a}, Section 4.1, we get
\begin{equation}
\rho(x,y) = Q(\zeta;\delta_1,\delta_2) := \frac{\ip{\cR(2\zeta+\delta_1+\delta_2)(\delta_1-\delta_2)}{\delta_1-\delta_2}}{\norm{2\zeta+\delta_1+\delta_2}^2 \ip{P_{J(2\zeta+\delta_1+\delta_2)}^{\perp}(\delta_1-\delta_2)}{\delta_1-\delta_2}} ~.
\end{equation}
Since $J(2\zeta+\delta_1+\delta_2) \in \mbox{ker}~\cR(2\zeta+\delta_1+\delta_2)$, we have
\begin{equation}
Q(\zeta;\delta_1,\delta_2) = \frac{\ip{\cR(2\zeta+\delta_1+\delta_2)P_{J(2\zeta+\delta_1+\delta_2)}^{\perp}(\delta_1-\delta_2)}{P_{J(2\zeta+\delta_1+\delta_2)}^{\perp}(\delta_1-\delta_2)}}{\norm{2\zeta+\delta_1+\delta_2}^2 \ip{P_{J(2\zeta+\delta_1+\delta_2)}^{\perp}(\delta_1-\delta_2)}{\delta_1-\delta_2}} ~.
\end{equation}
Now let $\delta = \delta_1+\delta_2$ and $\nu = \delta_1-\delta_2$. Note the set inclusion relation
\begin{equation*}
\begin{aligned}
& \left\{\delta_1,\delta_2 \in \RR^{2n}:~~ \norm{\delta}<\frac{r}{2},~ \norm{\nu}<\frac{r}{2},~ \nu \perp J(2\zeta+\delta) \right\} \\
\subset & \left\{\delta_1,\delta_2 \in \RR^{2n}:~~ \norm{\delta_1}<r,~ \norm{\delta_2}<r,~ \nu \perp J(2\zeta+\delta) \right\} \\
\subset & \left\{\delta_1,\delta_2 \in \RR^{2n}:~~ \norm{\delta}<2r,~ \norm{\nu}<2r,~ \nu \perp J(2\zeta+\delta) \right\} ~.
\end{aligned}
\end{equation*}
Thus we have
\begin{equation*}
\inf_{\substack{\norm{\delta}<2r \\ \norm{\nu}<2r \\ \nu \perp J(2\zeta+\delta)}} Q(\zeta;\delta_1,\delta_2) \quad \leq \inf_{\substack{\norm{\delta_1}<r \\ \norm{\delta_2}<r \\ \nu \perp J(2\zeta+\delta)}} Q(\zeta;\delta_1,\delta_2) \quad \leq \inf_{\substack{\norm{\delta}<r/2 \\ \norm{\nu}<r/2 \\ \nu \perp J(2\zeta+\delta)}} Q(\zeta;\delta_1,\delta_2) ~.
\end{equation*}
That is,
\begin{equation*}
\inf_{\norm{\delta}<2r} \frac{\l_{2n-1}(\cR(2\zeta+\delta))}{\norm{2\zeta+\delta}^2} \quad \leq \inf_{\substack{\norm{\delta_1}<r \\ \norm{\delta_2}<r \\ \nu \perp J(2\zeta+\delta)}} Q(\zeta;\delta_1,\delta_2) \quad \leq \inf_{\norm{\delta}<r/2} \frac{\l_{2n-1}(\cR(2\zeta+\delta))}{\norm{2\zeta+\delta}^2} ~.
\end{equation*}
Take $r \rightarrow 0$, by the continuity of eigenvalues with respect to the matrix entries, we have
\begin{equation*}
\l_{2n-1}(\cR(\zeta))/\norm{\zeta}^2 \leq a(z) \leq \l_{2n-1}(\cR(\zeta))/\norm{\zeta}^2 ~.
\end{equation*}
That is,
\begin{equation}
a(z) = \l_{2n-1}(\cR(\zeta))/\norm{\zeta}^2 ~.
\end{equation}

Now consider
\begin{equation}
\rho(x,z) = \frac{\norm{\beta(x)-\beta(z)}^2}{d_1(x,z)^2} ~.
\end{equation}
For simplicity write $\delta = \delta_1$. We can compute that
\begin{equation}
\rho(x,z) = Q(\zeta;\delta) = \frac{\ip{\cR(2\zeta+\delta)\delta}{\delta}}{\norm{2\zeta+\delta}^2 \ip{P_{J(2\zeta+\delta)}^{\perp} \delta}{\delta}} = \frac{\ip{\cR(2\zeta+\delta)P_{J(2\zeta+\delta)}^{\perp} \delta}{P_{J(2\zeta+\delta)}^{\perp} \delta}}{\norm{2\zeta+\delta}^2 \ip{P_{J(2\zeta+\delta)}^{\perp} \delta}{\delta}} ~.
\end{equation}
Note that
\begin{equation*}
\inf_{\substack{\norm{\delta}<r \\ \delta \perp J(2\zeta+\delta)}} Q(\zeta;\delta) \quad \geq \quad \inf_{\norm{\sigma}<r} \inf_{\substack{\norm{\delta}<r \\ \delta \perp J(2\zeta+\delta)}} Q(\zeta;\delta) \quad = \quad \inf_{\norm{\sigma}<r} \l_{2n-1} (\cR(2\zeta+\delta)) ~.
\end{equation*}
Take $r \rightarrow 0$ we have that
\begin{equation*}
\tilde{a}(z) \geq \l_{2n-1}(\cR(2\zeta))/\norm{2\zeta}^2 = \l_{2n-1}(\cR(\zeta))/\norm{\zeta}^2 ~.
\end{equation*}
On the other hand, take $\tilde{e}_{2n-1}$ to be a unit-norm eigenvector correspondent to $\l_{2n-1}(\cR(2\zeta))$. Then by the continuity of eigenvalues with respect to the matrix entries, for any $\epsilon > 0$, there exists $t>0$ so that $\delta = t \tilde{e}_{2n-1}$ satisfy
\begin{equation}
\frac{\ip{\cR(2\zeta+\delta)\delta}{\delta}}{\ip{P_{J(2\zeta+\delta)}^{\perp} \delta}{\delta}} \leq \l_{2n-1} (\cR(2\zeta)) + \epsilon
\end{equation}
and from there we have
\begin{equation}
\tilde{a}(z) \leq \l_{2n-1}(\cR(2\zeta))/\norm{2\zeta}^2 = \l_{2n-1}(\cR(\zeta))/\norm{\zeta}^2 ~.
\end{equation}
Therefore,
\begin{equation}
\tilde{a}(z) = \l_{2n-1}(\cR(\zeta))/\norm{\zeta}^2 ~.
\end{equation}

In a similar way (replacing infimum by supremum) we also get $b(z)$ and $\tilde{b}(z)$ as stated in the theorem.

Q.E.D.

\subsection{Proof of results in Section \ref{sec3}}\label{sec4_2}
We start by proving Proposition \ref{lem1}.

\noindent{\bf Proof of Proposition \ref{lem1}}

\begin{enumerate}

\item For $D_p$ obviously we have $D_p(\hat{x},\hat{y})\geq 0$ for any $\hat{x}$, $\hat{y} \in \hat{H}$ and $D_p(\hat{x},\hat{y})=0$ if and only if $\hat{x}=\hat{y}$. We also have $D_p(\hat{x},\hat{y})=D_p(\hat{y},\hat{x})$ since $\norm{x-ay}_{p}=\norm{y-a^{-1}x}_{p}$ for any $x$, $y \in H$, $|a|=1$. Moreoever, for any $x$, $y$, $z \in H$, if $D_p(\hat{x},\hat{y}) = \norm{x-ay}_{p}$, $D_p(\hat{y},\hat{z}) = \norm{z-by}$, then
\[
D_p(\hat{x},\hat{z}) \leq \norm{x-ab^{-1}z}_{p} = \norm{bx-az}_{p} \leq \norm{bx-aby}_{p} + \norm{aby-az}_{p} = D_p(\hat{x},\hat{y})+D_p(\hat{y},\hat{z}).
\]
Therefore $D_p$ is a metric.

$d_p$ is also a metric since $\norm{\cdot}_{p}$ in the definition of $d_p$ is the standard Schatten p-norm of a matrix.

\item For $p \leq q$, by H\"{o}lder's inequality we have for any $x=(x_1,x_2,...,x_n) \in H= \mathbb{C}^n$ that $\sum_{i=1}^{n}|x_i|^p \leq n^{(\frac{1}{p}-\frac{1}{q})}(\sum_{i=1}^{n}|x_i|^q)^{\frac{p}{q}}$. Thus $\norm{x}_{p} \leq n^{(\frac{1}{p}-\frac{1}{q})} \norm{x}_{q}$. Also since $\norm{\cdot}_{p}$ is homogeneous, if we assume $\norm{x}_{p}=1$ we have $\sum_{i=1}^{n}|x_i|^q \leq \sum_{i=1}^{n}|x_i|^p = 1$. Thus $\norm{x}_{q} \leq \norm{x}_{p}$. Therefore, we have $D_q(\hat{x},\hat{y})=\norm{x-a_1y}_{q} \geq n^{(\frac{1}{p}-\frac{1}{q})} \norm{x-a_1y}_{p} \geq n^{(\frac{1}{p}-\frac{1}{q})}D_p(\hat{x},\hat{y})$ and $D_p(\hat{x},\hat{y})=\norm{x-a_2y}_{p} \geq \norm{x-a_2y}_{q} \geq D_q(\hat{x},\hat{y})$ for some $a_1$, $a_2$ with magnitude $1$. Hence

\[D_q(\hat{x},\hat{y}) \leq D_p(\hat{x},\hat{y}) \leq n^{(\frac{1}{p}-\frac{1}{q})}D_q(\hat{x},\hat{y})\]

We see that $(D_p)_{1\leq p \leq \infty}$ are equivalent. The second part follows then immediately.

\item The proof is similar to (ii). Note that there are at most 2 $\sigma_i$'s that are nonzero, so we have $2^{(\frac{1}{p}-\frac{1}{q})}$ instead of $n^{(\frac{1}{p}-\frac{1}{q})}$.

\item To prove that $D_p$ and $d_q$ are equivalent, we need only to show that each open ball with respect to $D_p$ contains an open ball with respect to $d_p$, and vise versa. By (ii) and (iii), it is sufficient to consider the case when $p=q=2$.

First, we fix $x \in H = \mathbb{C}^n$, $r > 0$. Let $R = \min(1,\frac{r}{(2\norm{x}_{\infty}+1)n^2})$. Then for any $\hat{y}$ such that $D_2(\hat{x},\hat{y})<R$, we take $y$ such that $\norm{x-y}<R$, then $\forall 1 \leq i,j \leq n$, $|x_i\overline{x_j}-y_i\overline{y_j}| = |x_i(\overline{x_j}-\overline{y_j})+(x_i-y_i)\overline{y_j}| < |x_i|R+R(|x_i|+R) = R(2|x_i|+R) \leq R(2|x_i|+1) \leq \frac{r}{n^2}$. Hence $d_2(\hat{x},\hat{y}) = \norm{xx^*-yy^*}_{2} < n^2 \cdot \frac{r}{n^2} = r$.

On the other hand, we fix $x \in H = \mathbb{C}^n$, $R>0$. Let $r = R^2/\sqrt{2}$. Then for any $\hat{y}$ such that $d_2(\hat{x},\hat{y}) < r$, we have
\begin{equation}
 (d_2(\hat{x},\hat{y}))^2 = \norm{x}^4+\norm{y}^4-2|\langle x,y \rangle|^2 < r^2 = \frac{R^4}{2}
\end{equation}

But we also have
\begin{equation}
 (D_2(\hat{x},\hat{y}))^2 = \min_{|a|=1} \norm{x-ay}^2 = \norm{x-\frac{\langle x,y \rangle}{|\langle x,y \rangle|}y}^2 = \norm{x}^2+\norm{y}^2-2|\langle x,y \rangle|
\end{equation}

So
\begin{equation}
 (D_2(\hat{x},\hat{y}))^4 = \norm{x}^4 + \norm{y}^4 + 2\norm{x}^2 \norm{y}^2 - 4(\norm{x}^2+\norm{y}^2)|\langle x,y \rangle|+4|\langle x,y \rangle|^2
\end{equation}

Since $|\langle x,y \rangle| \leq \norm{x}\norm{y} \leq (\norm{x}^2+\norm{y}^2)/2$, we can easily check that $(D_2(\hat{x},\hat{y}))^4 \leq 2(d_2(\hat{x},\hat{y}))^2 < R^4$. Hence $D_2(\hat{x},\hat{y}) < R$.

Thus $D_2$ and $d_2$ are indeed equivalent metrics. Therefore $D_p$ and $d_q$ are equivalent. Also, the imbedding $i$ is not Lipschitz: if we take $x=(x_1,0,\ldots,0) \in \mathbb{C}^n$, then $D_2(\hat{x},0) = |x_1|$, $d_2(\hat{x},0) = |x_1|^2$.

\item First, for $p = 2$, for $\hat{x} \neq \hat{y}$ in $\hat{H}-\{0\}$, we compute the quotient
\begin{equation*}
\begin{aligned}
\rho(x,y) & = \frac{\norm{\kappa_{\alpha}(x)-\kappa_{\alpha}(y)}^2}{D_2(x,y)^2} \\
& = \frac{\norm{\norm{x}^{-1}xx^*-\norm{y}^{-1}yy^*}^2}{\norm{x}^2+\norm{y}^2-2\abs{\ip{x}{y}}} \\
& = \frac{\norm{xx^*}^2 \norm{y}^2 + \norm{x}^2 \norm{yy^*}^2 - 2\norm{x} \norm{y} \trace(xx^*yy^*)}{\norm{x}^4 \norm{y}^2 + \norm{x}^2 \norm{y}^4 - 2 \norm{x}^2 \norm{y}^2 \abs{x^* y}} \\
& = 1 + \frac{2\norm{x}\norm{y} \left( \norm{x}\norm{y}\abs{x^*y}-\trace(xx^*yy^*) \right) }{\norm{x}^4 \norm{y}^2 + \norm{x}^2 \norm{y}^4 - 2 \norm{x}^2 \norm{y}^2 \abs{x^* y}} \\
& = 1 + \frac{2 \left( \norm{x}\norm{y}\abs{x^*y}-\trace(xx^*yy^*) \right) }{\norm{x}^3 \norm{y} + \norm{x} \norm{y}^3 - 2 \norm{x} \norm{y} \abs{x^* y}}
\end{aligned}
\end{equation*}
where we used $\norm{xx^*} = \norm{x}^2$. For simplicity write $a = \norm{x}$, $b = \norm{y}$ and $t = \abs{\ip{x}{y}} \cdot (\norm{x}\norm{y})^{-1}$. We have $a > 0$, $b > 0$ and $0 \leq t \leq 1$.

Now we have
\begin{equation*}
\rho(x,y) = 1 + \frac{2(abt-abt^2)}{a^2+b^2-2abt}
\end{equation*}
Obviously, we have $\rho(x,y) \geq 1$. Now we prove that $\rho(x,y) \leq 2$. Note that
\begin{equation*}
1 + \frac{2(abt-abt^2)}{a^2+b^2-2abt} \leq 2 \Leftrightarrow a^2+b^2-4abt+2abt^2 \geq 0
\end{equation*}
But
\begin{equation*}
a^2+b^2-4abt+2abt^2 \geq 2ab-4abt+2abt^2 = 2ab(t-1)^2 \geq 0,
\end{equation*}
so we are done.
Note that take any $x$, $y$ with $\ip{x}{y} = 0$ we would have $\rho(x,y) = 1$. On the other hand, taking $\norm{x} = \norm{y}$ and let $t \rightarrow 1$ we see that $\rho(x,y) = 2-\epsilon$ is achievable for any small $\epsilon > 0$. Therefore the constants are optimal. The case where one of $x$ and $y$ is zero would not break the constraint of these two constants. Therefore after taking the square root we get lower Lipschitz constant $1$ and upper Lipschitz constant $\sqrt{2}$.

For other $p$, we use the results in (ii) and (iii) to get that the lower Lipschitz constant for $\kappa_{\alpha}$ is $\min(2^{\frac{1}{2}-\frac{1}{p}},n^{\frac{1}{p}-\frac{1}{2}})$ and the upper Lipschitz constant is $\sqrt{2}\max(n^{\frac{1}{2}-\frac{1}{p}},2^{\frac{1}{p}-\frac{1}{2}})$.

\item This follows directly from the construction of the map.

\item This follows directly from (v) and (vi).

Q.E.D.

\end{enumerate}

\vspace{5mm}

Next we prove Lemma \ref{prop:2}.

\noindent{\bf Proof of Lemma \ref{prop:2}}

(ii) follows directly from the expression of $\pi$. We prove (i) below.

Let $A$, $B \in \Sym(H)$ where $A=\sum_{k=1}^d \lambda_{m(k)}P_k$ and $B=\sum_{k'=1}^{d'} \mu_{m(k')}Q_{k'}$. We now show that
\begin{equation}
\label{eq:lipc}
 \norm{\pi(A)-\pi(B)}_{p} \leq (3+2^{1+\frac{1}{p}})\norm{A-B}_{p}
\end{equation}

Assume $\lambda_1-\lambda_2\leq \mu_1-\mu_2$. Otherwise switch the notations for $A$ and $B$. If $\mu_1-\mu_2=0$ then $\pi(A)=\pi(B)=0$ and
the inequality (\ref{eq:lipc}) is satisfied. Assume now $\mu_1-\mu_2>0$. Thus $Q_1$ is of rank 1 and therefore $\norm{Q_1}_p=1$ for all $p$.
First note that
\begin{equation}
\label{eq:break}
 \pi(A)-\pi(B) = (\lambda_1-\lambda_2)P_1 - (\mu_1-\mu_2)Q_1 = (\lambda_1-\lambda_2)(P_1-Q_1) + (\lambda_1-\mu_1-(\lambda_2-\mu_2))Q_1
\end{equation}

Here $\norm{P_1}_{\infty} = \norm{Q_1}_{\infty} = 1$. Therefore we have $\norm{P_1-Q_1}_{\infty} \leq 1$ since $P_1$, $Q_1 \geq 0$. From that we have $\norm{P_1-Q_1}_p \leq 2^{\frac{1}{p}}$.

Also, by Weyl's inequality we have $|\lambda_i-\mu_i| \leq \norm{A-B}_{\infty}$ for each $i$. Apply this to $i=1$, $2$ we get $|\lambda_1-\mu_1-(\lambda_2-\mu_2)| \leq |\lambda_1-\mu_1|+|\lambda_2-\mu_2| \leq 2\norm{A-B}_{\infty}$. Thus $|\lambda_1-\mu_1|+|\lambda_2-\mu_2| \leq 2\norm{A-B}_{\infty} \leq 2\norm{A-B}_p$.

Let $g:=\lambda_1-\lambda_2$, $\delta := \norm{A-B}_p$, then apply the above inequality to (\ref{eq:break}) we get
\begin{equation}
\label{eq:gp2d}
 \norm{\pi(A)-\pi(B)}_p \leq g\norm{P_1-Q_1}_p+2\delta \leq 2^{\frac{1}{p}}g+2\delta
\end{equation}

If $0 \leq g \leq (2+2^{-\frac{1}{p}})\delta$, then $\norm{\pi(A)-\pi(B)}_p \leq (2^{1+\frac{1}{p}}+3)\delta$ and we are done.

Now we consider the case where $g>(2+2^{-\frac{1}{p}})\delta$. Note that in this case we have $\delta < g/2$. Thus we have $|\lambda_1-\mu_1| < g/2$ and $|\lambda_2-\mu_2| < g/2$. That means $\mu_1 > (\lambda_1+\lambda_2)/2$ and $\mu_2 < (\lambda_1+\lambda_2)/2$. Therefore, we can use holomorphic functional calculus and put
\begin{equation}
 P_1 = -\frac{1}{2\pi i}\oint_\gamma R_Adz
\end{equation}
and
\begin{equation}
 Q_1 = -\frac{1}{2\pi i}\oint_\gamma R_Bdz
\end{equation}
where $R_A = (A-zI)^{-1}$, $R_B = (B-zI)^{-1}$, and $\gamma = \gamma(t)$ is the contour given in the picture below (note that $\gamma$ encloses $\mu_1$ but not $\mu_2$) and used also by \cite{ZB06}.

\begin{figure}[ht]
 \centering
 \includegraphics{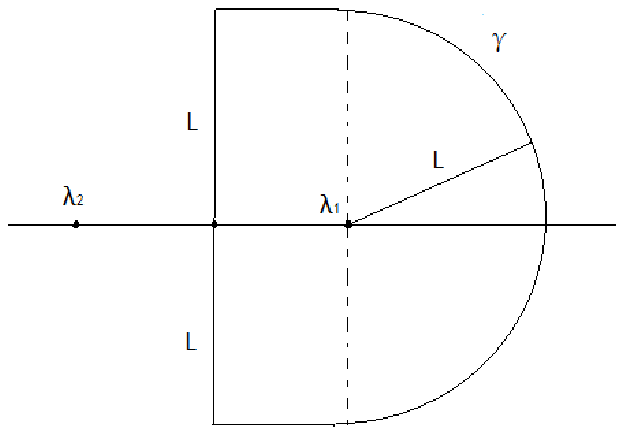}
\end{figure}

Therefore we have
\begin{equation}
\label{eq:pmq}
 \norm{P_1-Q_1}_p \leq \frac{1}{2\pi} \int_I \norm{(R_A-R_B)(\gamma(t))}_p |\gamma^{\prime} (t)| dt
\end{equation}

Now we have
\begin{equation}
\label{eq:ramrb}
 (R_A-R_B)(z) = R_A(z)-(I+R_A(z)(B-A))^{-1}R_A(z) = \sum_{n \geq 1} (-1)^n (R_A(z)(B-A))^n R_A(z)
\end{equation}
since for large $L$ we have $\norm{R_A(z)(B-A)}_{\infty} \leq \norm{R_A(z)}_{\infty}\norm{B-A}_p \leq \frac{\delta}{dist(z,\rho(A))} \leq \frac{2\delta}{g} < \frac{2}{2+2^{-\frac{1}{p}}} < 1$, where $\rho(A)$ denotes the spectrum of A.

Therefore we have
\begin{equation}
\begin{aligned}
\label{eq:nramrb}
 \norm{(R_A-R_B)(\gamma(t))}_p
 &\leq \sum_{n \geq 1} \norm{R_A(\gamma(t))}_{\infty}^{n+1} \norm{A-B}_p^n \\
 &=\frac{\norm{R_A(\gamma(t))}_{\infty}^2\norm{A-B}_p}{1-\norm{R_A(\gamma(t))}_{\infty}\norm{A-B}_p}
 <\frac{\norm{A-B}_p}{dist^2(\gamma(t),\rho(A))} \cdot (2^{1+\frac{1}{p}}+1)
\end{aligned}
\end{equation}
since $dist(\gamma(t),\rho(A)) \geq g/2$ for each t for large $L$. Here we used the fact that if we order the singular values of any matrix $X$ such that $\sigma_1(X) \geq \sigma_2(X) \geq \cdots$, then for any $i$ we have $\sigma_i(XY) \leq \sigma_1(X)\sigma_i(Y)$, and thus for two operators $X$, $Y \in \Sym(H)$, we have $\norm{XY}_p \leq \norm{X}_{\infty} \norm{Y}_p$.

Hence by (\ref{eq:pmq}) and (\ref{eq:nramrb}) we have
\begin{equation}
\label{eq:npmq}
 \norm{P_1-Q_1}_p \leq (2^{\frac{1}{p}}+2^{-1})\frac{\norm{A-B}_p}{\pi} \int_I \frac{1}{dist^2(\gamma(t),\rho(A))}|\gamma^{\prime}(t)|dt
\end{equation}

By evaluating the integral and letting $L$ approach infinity for the contour, we have as in  \cite{ZB06}
\begin{equation}
\label{eq:atan}
 \int_I \frac{1}{dist^2(\gamma(t),\rho(A))}|\gamma^{\prime}(t)|dt = 2\int_0^{\infty}\frac{1}{t^2+(\frac{g}{2})^2}dt = \left[\frac{4}{g}\arctan\left(\frac{2t}{g}\right)\right]_0^{\infty} = \frac{2\pi}{g}
\end{equation}

Hence
\begin{equation}
\label{eq:npmq2}
 \norm{P_1-Q_1}_p \leq (2^{\frac{1}{p}}+2^{-1})\frac{\norm{A-B}_p}{\pi} \cdot \frac{2\pi}{g} = (2^{1+\frac{1}{p}}+1)\frac{\delta}{g}
\end{equation}

Thus by the first inequality in (\ref{eq:gp2d}) and (\ref{eq:npmq2}) we have $\norm{\pi(A)-\pi(B)}_{p} \leq
(3+2^{1+\frac{1}{p}})\delta$.

We have proved that $\norm{\pi(A)-\pi(B)}_{p} \leq (3+2^{1+\frac{1}{p}})\norm{A-B}_{p}$.
That is to say, $\pi:(\Sym(H),\norm{\cdot}_p)\rightarrow(S^{1,0}(H),\norm{\cdot}_p)$
is Lipschitz continuous with Lipschitz constant less than or equal to $3+2^{1+\frac{1}{p}}$.

Q.E.D.
\vspace{5mm}

Now we are ready to prove Theorem \ref{maintheo}.

\noindent{\bf Proof of Theorem \ref{maintheo}}

The proof for $\alpha$ and $\beta$ are the same in essence. For simplicity we do it for $\beta$ first. 

We construct a map $\psi:(\R^m,\norm{\cdot}_p)\rightarrow (\hat{H},d_q)$ so that $\psi(\beta(x))=x$ for all $x\in\hat{H}$, and $\psi$ is Lipschitz continuous.
We prove the Lipschitz bound (\ref{eq:Lpq1}) which implies (\ref{eq:inv2}) for $p=2$ and $q=1$.

Set $M=\beta(\hat{H})\subset\R^m$. By hypothesis, there is a map $\tilde{\psi}_1:M\rightarrow \hat{H}$ that is Lipschitz continuous and satisfies
$\tilde{\psi}_1(\beta(x))=x$ for all $x\in\hat{H}$. Additionally, the Lipschitz bound between $(M,\norm{\cdot}_2)$ (that is, $M$ with Euclidian distance) and $(\hat{H},d_1)$
is given by $\frac{1}{\sqrt{a_0}}$.

First we change the metric on $\hat{H}$ from $d_1$ to $d_2$ and embed isometrically $\hat{H}$ into $\Sym(H)$ with Frobenius norm (i.e. Euclidian metric):
\begin{equation}
\label{eq:p1}
 (M,\norm{\cdot}_2) \stackrel{\tilde{\psi}_1}{\longrightarrow} (\hat{H},d_1) \stackrel{i_{1,2}}{\longrightarrow} (\hat{H},d_2) \stackrel{\kappa_\beta}{\longrightarrow}
(\Sym(H),\norm{\cdot}_{Fr})
\end{equation}
where $i_{1,2}(x)=x$ is the identity of $\hat{H}$ and $\kappa_\beta$ is the isometry (\ref{eq:kappa_b}) .
We obtain a map $\tilde{\psi}_2:(M,\norm{\cdot}_2)\rightarrow (\Sym(H),\norm{\cdot}_{Fr})$ of Lipschitz constant
$$\Lip(\tilde{\psi}_2)\leq \Lip(\tilde{\psi}_1)\Lip(i_{1,2})\Lip(\kappa_\beta) = \frac{1}{\sqrt{a_0}},$$
 where we used $\Lip(i_{1,2}) = L^d_{1,2,n}=1$ by (\ref{eq:L2}).

Kirszbraun Theorem \cite{WelWil75} extends isometrically $\tilde{\psi}_2$ from $M$ to the entire
$\R^m$ with Euclidian metric $\norm{\cdot}$. Thus we obtain a Lipschitz map $\psi_2:(\R^m,\norm{\cdot})\rightarrow (\Sym(H),\norm{\cdot}_{Fr})$ of Lipschitz constant
$\Lip(\psi_2)=\Lip(\tilde{\psi}_2)\leq \frac{1}{\sqrt{a_0}}$ so that $\psi_2(\beta(x))=\outp{x}{x}$ for all $x\in \hat{H}$.

The third step is to piece together $\psi_2$ with norm changing identities.

For $q\leq 2$ we consider the following maps:
\begin{equation}
\label{eq:p3}
(\R^m,\norm{\cdot}_p) \stackrel{j_{p,2}}{\longrightarrow}  (\R^m,\norm{\cdot}_2) \stackrel{{\psi}_2}{\longrightarrow} (\Sym(H),\norm{\cdot}_{Fr})
\stackrel{\pi}{\longrightarrow} (S^{1,0}(H),\norm{\cdot}_{Fr}) \stackrel{\kappa_\beta^{-1}}{\longrightarrow}
(\hat{H},d_{2}) \stackrel{i_{2,q}}{\longrightarrow} (\hat{H},d_q)
 \end{equation}
where $j_{p,2}$ and $i_{2,q}$ are identity maps on the respective spaces that change the metric.
The map $\psi$ claimed by Theorem \ref{maintheo} is obtained by composing:
\[ \psi:(\R^m,\norm{\cdot}_p) \rightarrow (\hat{H},d_q)~~,~~\psi = i_{2,q}\cdot \kappa_\beta^{-1} \cdot \pi \cdot \psi_2 \cdot j_{p,2} \]
Its Lipschitz constant is bounded by
\[ \Lip(\psi)_{p,q} \leq \Lip(j_{p,2}) \Lip(\psi_2)\Lip(\pi)\Lip(\kappa_\beta^{-1}) \Lip(i_{2,q}) \leq \max(1,m^{\frac{1}{2}-\frac{1}{p}}) \frac{1}{\sqrt{a_0}}\cdot (3+2\sqrt{2})\cdot 1
\cdot 2^{\frac{1}{q}-\frac{1}{2}} \]
Hence we obtained (\ref{eq:Lippq}). The other equation (\ref{eq:inv2}) follows for $p=2$ and $q=1$.

For $q>2$ we use:
\begin{equation}
\label{eq:p4}
(\R^m,\norm{\cdot}_p) \stackrel{j_{p,2}}{\longrightarrow}  (\R^m,\norm{\cdot}_2) \stackrel{{\psi}_2}{\longrightarrow} (\Sym(H),\norm{\cdot}_{Fr})
\stackrel{I_{2,q}}{\longrightarrow} (\Sym(H),\norm{\cdot}_{q})
\stackrel{\pi}{\longrightarrow} (S^{1,0}(H),\norm{\cdot}_{q}) \stackrel{\kappa_\beta^{-1}}{\longrightarrow}
(\hat{H},d_{q})
 \end{equation}
where $j_{p,2}$ and $I_{2,q}$ are identity maps on the respective spaces that change the metric.
The map $\psi$ claimed by Theorem \ref{maintheo} is obtained by composing:
\[ \psi:(\R^m,\norm{\cdot}_p) \rightarrow (\hat{H},d_q)~~,~~\psi = \kappa_\beta^{-1} \cdot \pi \cdot I_{2,q} \cdot \psi_2 \cdot j_{p,2} \]
Its Lipschitz constant is bounded by
\[ \Lip(\psi)_{p,q} \leq \Lip(j_{p,2}) \Lip(\psi_2)\Lip(I_{2,q})\Lip(\pi)\Lip(\kappa_\beta^{-1}) \leq \max(1,m^{\frac{1}{2}-\frac{1}{p}}) \frac{1}{\sqrt{a_0}}\cdot
 1 \cdot (3+2^{1+\frac{1}{q}}) \cdot 1 \]
Hence we obtained (\ref{eq:Lippq2}).

Replace $\beta$ by $\alpha$, $\psi$ by $\omega$, and $\kappa_\beta$ by $\kappa_\alpha$ in the proof above, using the Lipschitz constants for $\kappa_{\alpha}$ in Proposition \ref{lem1}, we obtain (\ref{eq:Lippq01}) and (\ref{eq:Lippq02}).

Q.E.D.

\section*{Acknowledgements}

The authors were supported in part by NSF grant DMS-1109498 and DMS-1413249. 
He also acknowledges fruitful discussions with
Krzysztof Nowak and Hugo Woerdeman (both from Drexel University) who pointed out several references, with Stanislav Minsker (Duke University)
for pointing out \cite{ZB06} and \cite{DavKah69}, and Vern Paulsen (University of Houston), Marcin Bownick 
(University of Oregon) and Friedrich Philipp (University of Berlin).

\end{document}